\documentclass[letterpaper, 11pt]{amsart}

\usepackage{amsmath,amsthm,amsfonts,amssymb,amscd}
\usepackage{bbm}
\usepackage{bm}
\usepackage{tikz}
\usepackage{tikz-cd}
\usepackage{appendix}
\usepackage{BOONDOX-calo}
\usepackage{standalone}
\usetikzlibrary{arrows,chains,matrix,positioning,scopes, cd}
\usepackage[letterpaper, left=3cm,right=3cm, top=3cm, bottom=3cm]{geometry}
\usepackage{adjustbox}
\usepackage{enumitem}
\usepackage[all]{xy}
\usepackage[
colorlinks=true, citecolor=blue, linkcolor=blue, urlcolor=red]{hyperref}
\usepackage{multicol}
\allowdisplaybreaks
\usepackage{mathrsfs}
\usepackage{here}
\usepackage{blkarray}
\usepackage{caption} 
\usepackage{footnote}
\usetikzlibrary{patterns}
\usepackage{accents}

\newcommand{\C}{{\mathbb C}}
\newcommand{\R}{{\mathbb R}}

\newcommand{\Z}{{\mathbb Z}}
\newcommand{\Q}{{\mathbb Q}}

\newcommand{\A}{{\mathbb A}}

\newcommand\dual{\raise0.9ex\hbox{$\scriptscriptstyle\vee$}}

\newcommand{\dbtilde}[1]{\accentset{\approx}{#1}}

\newcommand{\mathsym}[1]{{}}
\newcommand{\unicode}[1]{{}}

\newcommand{\ord}{\operatorname{ord}}
\newcommand{\res}{\operatorname{Res}}
\newcommand{\an}{\mathrm{an}}
\newcommand{\dR}{\mathrm{dR}}
\newcommand{\bG}{{\mathbf{G}}}
\newcommand{\wX}{{\widetilde{X}}}
\newcommand{\wA}{{\widetilde{A}}}
\newcommand{\wB}{{\widetilde{B}}}
\newcommand{\wD}{{\widetilde{\Delta}}}
\newcommand{\D}{{\Delta}}
\newcommand{\tw}{{w}}
\newcommand{\ow}{{\overline{\omega}}}

\theoremstyle{plain}
\newtheorem{thm}{Theorem} 
\newtheorem{prop}[thm]{Proposition}
\newtheorem{lemma}[thm]{Lemma}
\newtheorem{lem}[thm]{Lemma}

\numberwithin{thm}{subsection}
\numberwithin{equation}{section}

\newtheorem{manualtheoreminner}{Theorem}
\newenvironment{thm'}[1]{%
  \renewcommand\themanualtheoreminner{#1}%
  \manualtheoreminner
}{\endmanualtheoreminner}

\theoremstyle{definition}
\newtheorem{defn}[thm]{Definition}

\newtheorem{notation}[thm]{Notation}

\newtheorem{dfn}[thm]{Definition}

\theoremstyle{remark}
\newtheorem{rem}[thm]{Remark}
\newtheorem{example}[thm]{Example}
\newtheorem{rmk}[thm]{Remark}

\numberwithin{equation}{section}

\tikzset{>=stealth}

\makeatletter
\def\@seccntformat#1{%
  \protect\textup{\protect\@secnumfont
    \ifnum\pdfstrcmp{subsection}{#1}=0 \bfseries\fi
    \csname the#1\endcsname
    \protect\@secnumpunct
  }%
}  
\makeatother

\makeatletter
\@namedef{subjclassname@2020}{%
  $2020$ Mathematics Subject Classification}
\makeatother

\begin{document}

\title[Mixed motives and linear forms in the Catalan constant]{Mixed motives and linear forms in the Catalan constant}
\author{Payman Eskandari}
\address{Department of Mathematics and Statistics, University of Winnipeg, Winnipeg MB, Canada }
\email{p.eskandari@uwinnipeg.ca}
\author{V. Kumar Murty}
\address{Department of Mathematics, University of Toronto, Toronto ON, Canada }
\email{murty@math.toronto.edu}
\author{Yusuke Nemoto}
\address{Graduate School of Science, Chiba University, 
Yayoicho 1-33, Inage, Chiba, 263-8522 Japan}
\email{y-nemoto@waseda.jp}
\subjclass[2020]{11G99, 11M06, 14F42}
\begin{abstract}
We first give a geometric construction of a 2-dimensional mixed motive over $\Q$ with the Catalan constant $\bG=1-1/3^2+1/5^2-1/7^2+\cdots$ as a period. We then use this motive to obtain a supply of linear forms in 1 and $\bG$. We also explicitly compute the coefficients of 1 and $\bG$ in these linear forms.
\end{abstract}
\maketitle
\setcounter{tocdepth}{1}
\hypersetup{bookmarksdepth = 2}
\vspace{-.2in}
\tableofcontents

\section{Introduction}
\subsection{Motivation and context}
This paper fits in a program initiated by Francis Brown \cite{Br16}, which aims to use algebraic geometry to explain some of the known irrationality proofs and perhaps discover new ones. One standard approach to irrationality proofs is as follows: In order to prove that $\alpha\in \R$ is irrational, one looks for a sequence of nonzero linear forms
\[
\lambda_n:=a_n+b_n \, \alpha \quad \quad \quad (a_n,b_n\in\Q \, , \ n\in\Z_{>0})
\]
of sufficiently fast decay such that
\[
\lambda_n \cdot (\text{the lowest common denominator of $a_n$ and $b_n$})
\]
tends to zero as $n$ tends to infinity. If such linear forms can be found, irrationality of $\alpha$ will follow by an easy argument. A more sophisticated version of this approach due to Nesterenko gives a line of attack for finding lower bounds on the dimension of the $\Q$-span of a given set of numbers $\alpha_1,\ldots,\alpha_r$. As far as the authors know, all the known irrationality and $\Q$-linear independence proofs about odd Riemann zeta values (in particular, \cite{Apery}, \cite{Beukers}, \cite{Rivoal00}, \cite{BR}, \cite{Zudilin04}) are based on variations of this approach. 

When $\alpha$, or more generally, when $\alpha_1,\ldots,\alpha_r$ are periods of motives, Brown's philosophy proposes using algebraic geometry to produce linear forms in $\alpha_1,\ldots,\alpha_r$. This philosophy is supported by the evidence that almost all linear forms that appear in the known irrationality and $\Q$-linear independence proofs involving odd zeta values do indeed arise, although only in hindsight, from algebraic geometry. See \cite[\S 1]{Br16} and \cite[\S 1]{Dupont} for more details.

A difficulty in proving irrationality of odd zeta values is the fact that the linear forms that one constructs often involve too many unwanted terms. For example, one might get a linear form in 1 and all the Riemann zeta values $\zeta(n)$ where $n=2,\ldots, N$. Since one already knows $\Q$-linear independence of the numbers 1 and $\zeta(2m)\in \pi^{2m}\,\Q^\times$ ($m\in \Z_{>0}$) thanks to Lindemann, it becomes extremely difficult to obtain sufficiently strong bounds on the linear forms and the height of the coefficients to deduce anything nontrivial about odd zeta values. To rectify this, one wants to obtain linear forms with no unwanted terms. 

On the algebro-geometric side, the above challenge is related to finding a geometric construction of a motive {\it as small as possible} with the desired periods. Again, we refer to Brown's paper \cite{Br16} for more details. Here, the term {\it geometric} vaguely means to be built out of the relative cohomology groups of algebraic varieties using operations such as taking kernels and cokernels in an abelian category of motives.

For example, in order to prove irrationality of an odd zeta value $\zeta(2m+1)$ ($m\in\Z_{>0}$) using the above approach, one needs suitable linear forms in 1 and $\zeta(2m+1)$. It is known that there exists a mixed motive $Z_{2m+1}$ over $\Q$, unique up to isomorphism, that is an extension of $\Q(-(2m+1))$ by $\Q(0)$ and has a period matrix of the form
\[
\begin{pmatrix}
    1& \zeta(2m+1)\\
    0& (2\pi i)^{2m+1}
\end{pmatrix}.
\]
The proof of this fact (which combines deep works of Beilinson \cite{Beilinson}, Borel \cite{Borel} and others; see \cite{DG}) does not give a geometric construction of this motive. Brown's philosophy proposes that a geometric construction of $Z_{2m+1}$ could lead to useful linear forms in 1 and $\zeta(2m+1)$. Interestingly, the only odd zeta value for which such a geometric construction is known is $\zeta(3)$. In that case, two different constructions are known, one due to Brown \cite{Br16} and another due to Dupont \cite{Dupont}. In fact, Dupont's work gives, for every integer $m\geq 1$, a geometric construction of a mixed motive in
\[
\textup{Ext}^1\bigm(\Q(-3)\oplus \Q(-5)\oplus\cdots\oplus \Q(-(2m+1)) , \, \Q(0)\bigm)
\]
whose periods include $\zeta(3), \zeta(5), \ldots, \zeta(2m+1)$. The Ball--Rivoal (\cite{Rivoal00}, \cite{BR}) linear forms in odd zeta values arise, in hindsight, from this motive. A functional version of \cite{Dupont} due to Dupont and Fres\'{a}n \cite{DF} gives a geometric construction of the polylogarithm motive.

In this paper we carry out an important part of Brown's program for the Catalan constant
\[
\bG  \ := \ 1-\frac{1}{3^2}+\frac{1}{5^2}-\frac{1}{7^2}+\cdots =L(2,\large{\chi}_{_{-4}}),
\]
where $L(s,\large{\chi}_{_{-4}})$ is the Dirichlet $L$-function for the nontrivial character mod 4. We give a geometric construction of a mixed motive $C$ over $\Q$ with a period matrix
\begin{equation}\label{p1eq11}
\begin{pmatrix}
1 & \bG\\
0 & \pi^2 i
\end{pmatrix}.
\end{equation}
The motive is not a mixed Tate motive over $\Q$, but over the field $\Q(i)$ it becomes an extension of $\Q(-2)$ by $\Q(0)$.\footnote{Recall that $\mathrm{Ext}^1(\Q(-2),\Q(0))$ should vanish for mixed motives over $\Q$, but not for mixed motives over $\Q(i)$.} We use the motive $C$ to produce a large supply of linear forms in 1 and $\bG$:
\begin{thm}\label{thm: linear forms}
Let $t$ be an integer $\geq 0$. Suppose $F=F(x,y)\in \Q[x,y]$ satisfies conditions (i) and (ii) below:
\begin{enumerate}
    \item $F$ is fixed by the automorphism $\sigma$ of $\Q[x,y]$ that sends $(x,y)\mapsto (-y,x)$.
    \item $x^{2\lceil t/2\rceil}y^{2\lceil t/2\rceil}$ (where $\lceil  \ \rceil$ is the ceiling function) divides $F$ in $\Q[x,y]$. 
\end{enumerate}
\noindent Then the integral of $\displaystyle{\frac{F dxdy}{(1-x^2-y^2)^{t+1}}}$ over the simplex
\[
\Delta:=\{(x,y)\in \R^2: x,y, 1-x-y \geq 0\}
\]
converges to a $\Q$-linear combination of $1$ and $\bG$.
\end{thm}

With potential applications to rational approximations and irrationality in mind, we make the proof of this result completely constructive, in that given $F$ and $t$ as above, we compute both coefficients of 1 and $\bG$ in the associated linear form explicitly. Denoting the least common multiple of $1,\ldots,m$ by $L_m$, in the case $t=0$, {for every polynomial $F$ as above and with coefficients in $\Z$, one has an upper bound of
\[
2^{\deg(F)+2} L_{\deg(F)}\,L_{\deg(F)/2} 
\]
for the lowest common denominator of the coefficients of 1 and $\bG$ in terms of the (total) degree of $F$. This does not seem to be good enough to prove irrationality of $\bG$ (see Remark \ref{rem: chebychev}), as one might expect since this bound is uniform for all $F\in \Z[x,y]$ of a given degree, and also one typically makes use of poles of higher order to obtain good rational approximations. The question remains if clever choices of $(F,t)$ can bring one any closer to an irrationality proof. We hope to pursue this line of inquiry in the future. It will also be interesting to compare the linear forms given by Theorem \ref{thm: linear forms} with those obtained by Rivoal and Zudilin (\cite[\S 9]{RZ} and \cite{Rivoal06}) and Nesterenko \cite{Nesterenko}.

Here are some examples of linear forms of Theorem \ref{thm: linear forms}:
\begin{align*}
&\int_\D \frac{x^2y^2}{1-x^2-y^2}dxdy \ = \ -\frac{5}{48} \, + \frac{1}{8} \, \bG \\
&\int_\D \frac{x^4y^4}{1-x^2-y^2}dxdy \ = \ -\frac{569}{26880}+\frac{3}{128} \bG \\
&\int_\D \frac{x^4y^4}{(1-x^2-y^2)^3}dxdy \ = \ -\frac{49}{384} + \frac{9}{64} \, \bG.
\end{align*}
The $\Z$-linear form obtained by clearing the denominators of the coefficients in each form (which is the quantity that we need to make arbitrarily close to zero by appropriate choices of $F$ and $t$ to deduce irrationality of $\bG$) is roughly $0.5$, $8.06$ and $0.46$, in that order. For $F=x^4y^4$, the choice $t=2$ (leading to the last form above) is the best possible choice between $0\leq t\leq 4$ ( = the range allowed for $t$ to assure convergence).

\subsection{Contents of the paper}
Our departure point is the integral expression\footnote{For a reference for this expression, see \cite{JL}, equation (20), where a variant of this formula is established. The variant used here is easily seen to be equivalent to the formula established in op. cit.}
\[
\bG \ = \int_\D \frac{dxdy}{1-x^2-y^2},
\]
where $\D$ is the simplex defined in Theorem \ref{thm: linear forms}. Using standard ideas, we use this integral expression to define a mixed motive\footnote{Say, in the setting of mixed absolute Hodge motives (defined in terms of mixed realizations, see \cite{HM17}, for instance) or Nori motives. See \S \ref{sec: settings} for more details.} over $\Q$ that has $\bG$ as a period. Section \ref{construction} gives the construction of this motive, denoted by $M$, as a relative cohomology group. Section \ref{sec: wt filtration} studies the weight filtration on $M$. It turns out that $M$ is 5-dimensional (i.e., its (say) Betti realization is 5-dimensional), and its associated graded with respect to the weight filtration is
\[
\operatorname{Gr}^WM \ \simeq \ \Q(0)\oplus \Q(-1)^{\oplus 3}\oplus H^1(A)(-1), 
\]
where $A$ is the circle $x^2+y^2=1$ inside the affine plane over $\Q$. Over $K=\Q(i)$, $A$ is isomorphic to $\mathbb{G}_m$ and hence the base change of $M$ to a motive over $K$ is a mixed Tate motive. In particular, $M$ is a mixed Artin--Tate motive over $\Q$ (by definition of such motives).

Section \ref{sec: Betti-dR realization} is devoted to the Betti and de Rham realizations of $M$. We give explicit bases for each of these and compute the periods of $M$, which turn out to be $\Q$-linear combinations of 1, $\log2$, $\bG$, $2\pi i$, and $\pi^2i$.

Section \ref{Catalanmotive} gives the definition of what we call the {\it motive of the Catalan constant}, denoted by $C$, which is a submotive of $M$ (over $\Q$) with period matrix \eqref{p1eq11} with respect to our bases. The key observation here is that one can use a particular symmetry (namely, $\sigma$ of Theorem \ref{thm: linear forms}) to precisely extract the desired submotive from $M$. The symmetry $\sigma$ does not act on the motive $M$ itself, but rather, acts on a quotient of $M$. A similar situation also occurs in \cite{Dupont}.

Sections \ref{sec: lin forms I}-\ref{sec: lin forms III} are devoted to using the motive $C$ to obtain linear forms in 1 and $\bG$. The goal is to prove Theorem \ref{thm: linear forms}, and compute the coefficients of 1 and $\bG$ in the linear forms (given $F$ and $t$). The idea is to find differential forms that give rise to elements of the de Rham realization of $C$. This is done in \S \ref{sec: lin forms I} (Proposition \ref{prop: proof of lin form thm, t=0 case}), leading to a non-constructive motivic proof of Theorem \ref{thm: linear forms}. In fact, we prove a more intrinsic version of Theorem \ref{thm: linear forms} (see Theorem \ref{thm: linear forms, ver2}). This second version seems at first to be stronger, but turns out to be equivalent to the version stated earlier, even for the purposes of calculation of the coefficients of 1 and $\bG$ (see \S \ref{sec:integrability}, in particular, the comments after Theorem \ref{thm: linear forms, ver2}, including Remark \ref{rem: ver 1 implies ver 2}). In \S \ref{sec: lin forms II} and \S \ref{sec: lin forms III} we make Theorem \ref{thm: linear forms} (and Theorem \ref{thm: linear forms, ver2}) constructive by calculating the coefficients of 1 and $\bG$. Section \ref{sec: lin forms II} handles the logarithmic case ($t=0$). The computation of the coefficient of $\bG$ in this case amounts to an easy residue calculation (see Proposition \ref{prop: computation of b_F}). The computation of the coefficient of 1 is considerably more involved. Nonetheless, we calculate this coefficient in \S \ref{sec: coefficient of eta, log case}. The approach is completely explicit. Combining formulas \eqref{p2eq9}, \eqref{p2eq6}, and Lemmas \ref{lem: 1st type integrals, log case} and \ref{lem: 2nd type integrals, log case} one can easily compute the coefficient of 1 in the linear forms in the logarithmic case. One also obtains bounds on the size of the denominators (see Theorem \ref{thm: lin forms thm, t=0, explicit version}). Finally, in \S \ref{sec: lin forms III} we make the case $t>0$ of the theorem constructive, by using an elementary procedure to reduce the order of the pole first.

\section{The motivic setting and notation}\label{sec: settings}
Let $F$ be a subfield of $\C$. By a variety over $F$ (or simply, a variety, when $F$ is understood from the context) we mean a quasi-projective reduced scheme of finite type over $F$. Throughout the paper, we work in a suitable category $\mathbf{MM}(F)$ of mixed motives over $F$ (specifically, we shall need $F=\Q$ and occasionally $F=\Q(i)$). To be concrete, the reader can take this to be the Tannakian category of (non-effective) Nori motives over $F$ with rational coefficients \cite{HM17}, or the category of mixed absolute Hodge motives over $F$ as in \cite[\S 6.3]{HM17} and \cite{Hub95} (note that these versions of the theory include motives of non-smooth varieties as well as motives for relative (co)homology). What we precisely require from $\mathbf{MM}(F)$ is the following: 
\begin{enumerate}[label=(\arabic*)]
\item $\mathbf{MM}(F)$ is a rigid abelian tensor category in the sense of \cite{DM} with $\mathrm{End}(\mathbbm{1})=\Q$, where $\mathbbm{1}$ is the unit object for the tensor product. (Note, importantly, that we work with motives with coefficients in $\Q$.)
\item There is a weight filtration $W_\bullet$ on objects of $\mathbf{MM}(F)$. This means that $W_\bullet$ is indexed by $\Z$, increasing, finite on every object, functorial, exact, and compatible with the tensor structure. (Functoriality and exactness amount to strict compatibility of morphisms with the weight filtration.)
\item There exist functorial dual objects $H^n(X)$ and $H_n(X)$ in $\mathbf{MM}(F)$ for the (co)homology of varieties (over $F$). There are also functorial dual objects $H^n(X,Z)$ and $H_n(X,Z)$ for the relative (co)homology of pairs $(X,Z)$ of a variety $X$ and a closed subvariety $Z\subset X$. When the second entry of the pair is the empty subvariety, the relative (co)homology recovers the absolute (co)homology. We have $H_0(\operatorname{Spec} F)=\mathbbm{1}$. 
\item There exist Betti, de Rham and Hodge realization functors with values respectively in the categories of vector spaces over $\Q$ and $F$, and the category of rational mixed Hodge structures. These functors are exact, faithful, $\Q$-linear, and compatible with the tensor structures (in particular, they are conservative). Thus the Betti realization makes $\mathbf{MM}(F)$ a neutral Tannakian category over $\Q$ in the sense of \cite{DM}. The realizations of the (co)homology motives are exactly those expected: the Betti realizations of $H^n(X,Z)$ and $H_n(X,Z)$ are the (degree $n$) rational singular cohomology and homology of the analytification of $(X,Z)$; the de Rham realization of $H^n(X,Z)$ is the relative algebraic de Rham cohomology of $(X,Z)$; and the Hodge realization is Deligne's mixed Hodge structure on the rational singular (co)homology. The Hodge realization functor translates the motivic weight to the Hodge-theoretic weight. In particular, the Lefschetz motive $\Q(-1):=H^2(\mathbb{P}^1)$ is pure of weight 2 and 1-dimensional in each of its realizations.
\item The following long exact sequences exist in $\mathbf{MM}(F)$:
\begin{itemize}
    \item[-] the long exact sequences for the relative cohomology and homology of a pair $(X,Z)$ as above, related to each other by duality;
    \item[-] the Mayer--Vietoris long exact sequence for an open covering $X=U\cup V$ of a variety $X$;
    \item[-] the Gysin long exact sequence for any triple $(X,Z,A)$ of a smooth variety $X$, a closed subvariety $Z\subset X$, and a smooth closed subvariety $A\subset X$ of pure codimension.
\end{itemize}
The Betti and de Rham realizations of these sequences are the corresponding sequences in the relevant theory.
\item There is a distinguished isomorphism called the period isomorphism between the complexifications of the Betti and de Rham realization functors, extending the usual Betti--de Rham comparison isomorphism for the motives $H^n(X)$ of smooth varieties to all objects of $\mathbf{MM}(F)$. In the case of the motives of the form $H^n(X,Z)$, the associated pairing (i.e., the period pairing) between the Betti realization of $H_n(X,Z)$ and the de Rham realization of $H^n(X,Z)$ is as in \cite[\S 5.3 -- \S 5.5]{HM17}.
\end{enumerate}

An excellent reference for almost all of the above in the settings of absolute Hodge and Nori is \cite{HM17}. In particular, see \S 3.2 and \S 3.3 therein for a detailed exposition of the relative algebraic de Rham cohomology. Existence of the long exact sequences of item (5) in the setting of Nori motives (which by Proposition 10.1.2 of \cite{HM17} and its proof also implies the existence of these exact sequences in the absolute Hodge setting) is probably well-known to some experts. For a reference, we refer the reader to \cite{FJ}, where the statements are proved in the more general setting of the authors' exponential motives (see op. cit., \S 4.2.4 for the relative (co)homology sequence, Corollary 4.6.5 for the Mayer--Vietoris sequence, and Theorem 4.8.3 for the Gysin sequence). By \S 5.1 of \cite{FJ}, taking the potential function to be zero, one recovers the desired statements for Nori motives from those for exponential motives.

We end this discussion with a few words about our notation. The complex analytic space corresponding to a variety $X$ over a subfield of $\C$ is denoted by $X^{\an}$. The Betti (resp. de Rham) realization of a motive $N$ is denoted by $N_{\textup{sing}}$ (resp. $N_\textrm{dR}$). Thus for example, the notation $H_{n,\textup{sing}}(X,Z)$ means the Betti realization of $H_n(X,Z)$, which is the relative singular homology $H_n(X^{\an}, Z^{\an};\Q)$, and the notation $H^n_{\mathrm{dR}}(X,Z)$ means the de Rham realization of the motive $H^n(X,Z)$, which is the degree $n$ algebraic de Rham cohomology of the pair $(X,Z)$. Similarly, the Betti (resp. de Rham) realization of a morphism $f$ of motives is denoted by $f_{\textup{sing}}$ (resp. $f_{\textup{dR}}$). Following the conventions of \cite{Deligne89}, the Betti (resp. de Rham) realization of $\Q(-1)$ over $\Q$ is $(2\pi i)^{-1}\Q$ (resp. $\Q$), with the comparison isomorphism given by the inclusions in $\C$.

Finally, we use the notation $V^\vee$ (resp. $f^\vee$) for the dual of an object $V$ (resp. morphism $f$) in a Tannakian category.

\section{Construction of the motive $M$} \label{construction}
\numberwithin{thm}{section}
In this section, we define a motive $M$ over $\Q$ that has $\bG$ as a period. This will be our first approximation of the motive of the Catalan constant, which will later be defined as a 2-dimensional subobject of $M$. The ideas used in the construction of $M$ go back to Terasoma \cite{Terasoma}. An excellent exposition of these ideas can be found in \cite[\S 2.5.1]{BF}, where Burgos Gil and Fres\'{a}n discuss the construction of a motive with $\zeta(2)$ as a period in details. This section will closely follow Burgos Gil and Fres\'{a}n.
\medskip\par 
Let us start by fixing some notation. From this point on in the paper, let $X:=\A^2=\operatorname{Spec}{\Q}[x, y]$ be the affine plane over $\Q$. Let $A \subset X$ be the unit circle around the origin, i.e., $$A:= \operatorname{Spec}\Q[x, y]/(1-x^2-y^2). $$
Let $\ell_1, \ell_2$ and $\ell_3$ be the lines in $X$ defined by 
\begin{align*}
&\ell_1=\{y=0\}, \quad \ell_2=\{x=0\}, \quad \ell_3=\{x+y=1\}, 
\end{align*}
and set $B=\ell_1 \cup \ell_2 \cup \ell_3$. We denote the point $(1,0)$ (resp. $(0,1)$) of $X$ by $P$ (resp. $Q$). We shall set (as we already did in the introduction)
$$\D:=\{(x, y)\in\R^2 : x, y ,1-x-y \geq 0\} \subset X^{\an}.$$ 
This is a topological $2$-chain in $X^{\an}$ with boundary in $B^\an$.
(See Figure \ref{fig1}.)

\begin{figure}[h]
\begin{tikzpicture}
\draw(0,0) circle (1.5cm);
\draw[->] (-2,0) -- (2, 0);
\draw[->] (0,-2)-- (0,2);
\draw[-] (-0.5, 2) -- (2, -0.5);
 \coordinate[] (A) at (0,0); 
 \coordinate[] (B) at (1.5,0); 
 \coordinate[] (C) at (0,1.5); 
\draw[pattern=north east lines] (0,0) -- (1.5,0) -- (0,1.5) -- cycle;
\coordinate[label=left:$\ell_1$] (l_1) at (-2.1,0);
 \coordinate[label=left:$\ell_2$] (l_2) at (0, -2.1);
  \coordinate[label=left:$\ell_3$] (l_3) at (2.1, -0.6);
    \coordinate[label=left:$P$] (P) at (2, 0.2);
     \coordinate[label=left:$Q$] (Q) at (0.5, 1.7);
     \foreach \P in {B,C} \fill[] (\P) circle[radius=0.06];
     \fill[white] (0.3,0.3) rectangle (0.7,0.7);
      \coordinate[label=left:$\Delta$] (D) at (0.8, 0.5);
\end{tikzpicture}
\hspace{7mm}
\vspace{-3mm}
\caption{}
\label{fig1}
\end{figure}

The starting point of the definition of $M$ is the integral expression
\[
\bG = \int_\Delta \frac{dxdy}{1-x^2-y^2}. 
\]
The differential form $dxdy/(1-x^2-y^2)$ is regular on $X\setminus A$. One may be tempted to immediately say that $\bG$ is a period of the motive
\[
H^2(X\setminus A, B\setminus A), 
\]
but $\Delta$ intersects $A$ at $P$ and $Q$. This is a standard issue (e.g., it is similar to the situation for multiple zeta values), and the remedy is straightforward: one considers the blow-up at the points $P$ and $Q$, as follows.
\medskip\par 
Let
$$\pi \colon \widetilde{X} \to X$$
 be the blow-up of $X$ at $P$ and $Q$. Explicitly, throughout the paper, we shall take $\wX$ to be the subvariety of 
 \begin{equation}\label{p1eq12}
 X \times \mathbb{P}^1 \times \mathbb{P}^1=\A^2 \times \mathbb{P}^1 \times \mathbb{P}^1=\{(x, y)\} \times \{[\alpha_1 : \beta_1]\} \times \{[\alpha_2 : \beta_2]\}
 \end{equation}
 defined by 
$$\left\{
\begin{array}{ll}
\alpha_1 \, (x-1)=\beta_1 \, y   \\
\alpha_2 \,x= \beta_2 \,(y-1). 
\end{array}
\right.
$$
The blow-up map $\pi \colon \wX \to X$ is the projection to the first component. We denote the special divisors above $P$ and $Q$ by $E_P$ and $E_Q$: 
\begin{align*}
&E_P =\pi^{-1}(P)=\pi^{-1}((1, 0))=\{(1, 0)\} \times \mathbb{P}^1 \times \{[1 : -1]\}  \\ 
&E_Q =\pi^{-1}(Q)=\pi^{-1}((0, 1))=\{(0, 1)\} \times \{[1 : -1]\} \times \mathbb{P}^1. 
\end{align*}
Set $E=E_P \cup E_Q$. We write $\widetilde{A}$ (resp. $\widetilde{B}$) for the strict transform of $A$ (resp. $B$). A straightforward computation shows
\begin{align}
\begin{split} \label{eq-tildeA cap E}
\wA\cap E_P & = \{P\} \times \{[1: 0]\} \times \{[1: -1]\}\\
\wA\cap E_Q & = \{Q\} \times \{[1: -1]\} \times \{[0: 1]\}.
\end{split}
\end{align}

\begin{defn}\label{defn of M}
Let $M$ be the motive
$$H^2(\wX \setminus \wA, (\wB \cup E) \setminus \wA)$$
in $\mathbf{MM}(\Q)$.
\end{defn}
\begin{rem}\label{rem: Xtilde-Atilde is affine}
The variety $\wX\setminus\wA$ is smooth and affine. Indeed, smoothness is clear and the reader can check\footnote{We include a proof sketch for convenience. The $f_j$ generate the ideal of $\{P,Q\}$ in $\Q[x,y]$. The blow-up $\wX$ has an affine open cover $X=\cup_{j=1}^4 V_j$ where $V_j=\operatorname{Spec}\Q[x,y,f_1/f_j, \ldots, f_4/f_j]$. Setting $g:=1-x^2-y^2$, the strict transform $\wA$ is defined in $V_j$ by $g/f_j = 0$ (note that $g=-f_1-f_2+f_3-f_4$ so $g/f_j \in \mathcal{O}(V_j)$), so that
\[
V_j\setminus\wA = \operatorname{Spec} \Q[x,y,f_1/f_j, \ldots, f_4/f_j, f_j/g] = \operatorname{Spec}\Q[x,y,f_1/g, \ldots, f_4/g, g/f_j]
\]
(recall $g/f_j \in \mathcal{O}(V_j)$). Thus $V_j\setminus\wA$ is simply the principal affine open set $D(f_j/g)$ of $\operatorname{Spec}\Q[x,y,f_1/g, \ldots, f_4/g]$ given by $f_j/g\neq 0$. One checks that the inclusions of the $V_j\setminus\wA$ in $\operatorname{Spec}\Q[x,y,f_1/g, \ldots, f_4/g]$ glue together to give an isomorphism between $\wX\setminus\wA$ and the union of the $D(f_j/g)$ in $\operatorname{Spec}\Q[x,y,f_1/g, \ldots, f_4/g]$, which is all of the latter because $1\in (f_1/g,\ldots, f_4/g)$.} that with $f_1=x(x-1)$, $f_2=xy$, $f_3=(x-1)(y-1)$, $f_4=y(y-1)$, we have
\[
\wX\setminus\wA \ \simeq \ \operatorname{Spec}\Q[x,y, \frac{f_1}{1-x^2-y^2}, \frac{f_2}{1-x^2-y^2},\frac{f_3}{1-x^2-y^2},\frac{f_4}{1-x^2-y^2}].
\]
\end{rem}

In the remainder of this section we carefully check that $\bG$ is a period of $M$. In the process, we will also fix some further notation.

First, we define an element of the Betti homology
\begin{equation}\label{p1eq1}
H_{2, \textup{sing}}(\wX \setminus \wA, (\wB \cup E) \setminus \wA) = M^\vee_{\textup{sing}}.
\end{equation}
Let $\wD$ be the closure (in the analytic topology) of $\pi^{-1}(\D \setminus \{P, Q\})$ in $\wX^\an$ (see Figure \ref{fig: wD}). This is a topological $2$-chain in $\wX^\an$ with $\pi(\wD)=\D$.

\begin{figure}
    \centering
    \includegraphics[width=0.5\linewidth]{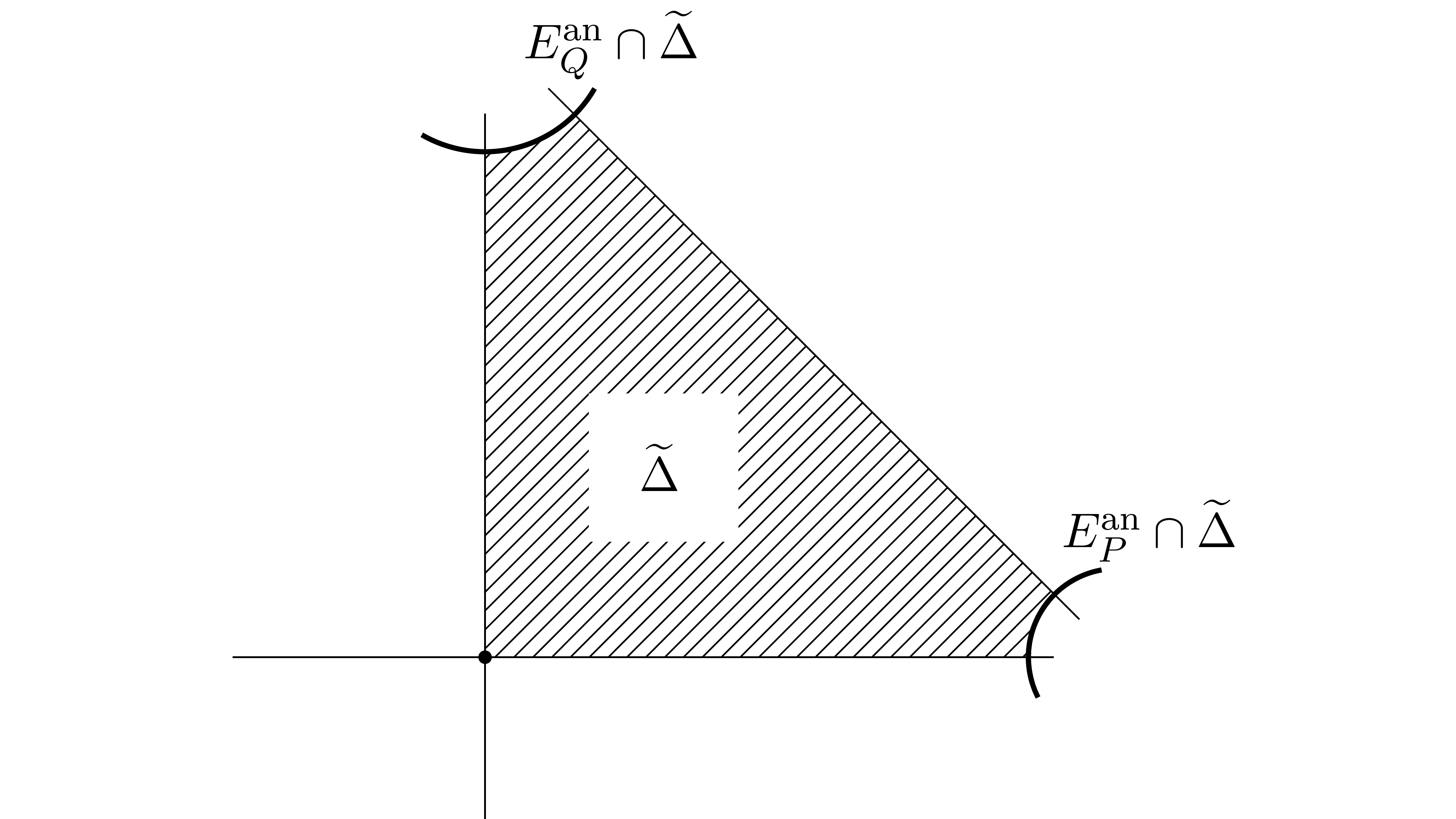}
    \vspace{-3mm}
    \caption{The cycle $\wD\subset \wX^\an$} 
    \label{fig: wD}
\end{figure}

\begin{lem} \label{Betti}
The topological chain $\wD$ is contained in $\wX^\an\setminus\wA^\an$ and has boundary in $\wB^\an\cup E^\an$.
\end{lem}

\begin{proof}
It is clear that away from $E_P^\an$ and $E_Q^\an$, the topological chain $\wD$ does not intersect $\wA^\an$. We have
\begin{align*}
&\widetilde{\Delta} \cap E_P^\an= \{P\} \times \{[m: 1] \mid -1 \leq m \leq 0\} \times \{[1: -1]\} \\ 
&\widetilde{\Delta} \cap E_Q^\an= \{Q\} \times \{[1: -1]\} \times \{[1: m] \mid -1 \leq m \leq 0\}.
\end{align*}
In particular, on recalling \eqref{eq-tildeA cap E}, $\wD$ does not intersect with $\wA^\an$ along the special divisors either, so that $\wD$ is contained in $\wX^\an\setminus\wA^\an$.
Moreover, we have
\begin{align} \label{boundary}
\partial \widetilde{\Delta} \subset \widetilde{\ell_1}^\an \cup \widetilde{\ell_2}^\an \cup \widetilde{\ell_3}^\an \cup (\widetilde{\Delta} \cap E_P^\an) \cup (\widetilde{\Delta} \cap E_Q^\an) \subset \wB^\an \cup E^\an, 
\end{align}
where $\partial$ refers to the boundary and $\widetilde{\ell_i}$ is the strict transform of $\ell_i$ in $\wX$.
\end{proof}
In view of the previous lemma, $\wD$ defines an element of the Betti homology
\eqref{p1eq1}, which with abuse of notation we also denote by $\wD$.

Next, we check that the differential form 
\begin{equation}\label{p1eq2}
\omega_0:=\pi^*\left(\dfrac{dxdy}{1-x^2-y^2} \right)=\dfrac{dxdy}{1-x^2-y^2}
\end{equation}
defines an element of the de Rham realization of $M$.

\begin{lem} \label{deRham}
The differential form $\omega_0$ on $\wX$ defined above is regular on $\wX \setminus \wA$, with a logarithmic singularity (i.e., pole of order $1$) along $\wA$.
\end{lem}

\begin{proof}
This is checked by a computation in charts. Consider the open covering 
$$\wX=\bigcup_{j=1}^4 U_j$$
of $\wX$, where the $U_j$ are coordinate charts as follows\footnote{This covering of $\wX$ is more convenient for some of the later arguments in the paper than the covering of $\wX$ by the open sets $\operatorname{Spec}\Q[x,y,f_1/f_j, \ldots, f_4/f_j]$ ($1\leq j\leq 4$), with the $f_j$ as in Remark \ref{rem: Xtilde-Atilde is affine} (forming a generating set of the ideal of $\{P,Q\}$).}:
\begin{align*}
&U_1= \wX\setminus (\{\alpha_1 =0\}\cup E_Q) \ \text{with affine coordinates} \ (y,v=\beta_1/\alpha_1) \\
&U_2= \wX\setminus (\{\beta_1 =0\}\cup E_Q) \ \text{with affine coordinates} \ (x, u=\alpha_1/\beta_1) \\
&U_3= \wX\setminus (E_P\cup \{\alpha_2 =0\}) \ \text{with affine coordinates} \ (y,u=\beta_2/\alpha_2)\\
&U_4= \wX\setminus (E_P\cup \{\beta_2 =0\}) \ \text{with affine coordinates} \ (x, v=\alpha_2/\beta_2).
\end{align*}
On $U_1$, $x=1+vy$, 
\[
1-x^2-y^2 = y(-v^2y -2v -y),
\]
and $\wA$ is given by $-v^2y -2v -y=0$. The differential form $\omega_0$ is written as
$$\dfrac{dxdy}{1-x^2-y^2}=\dfrac{dvdy}{-v^2 y-2v-y}, $$
which has a logarithmic singularity along $\wA$ and is regular outside $\wA$. On $U_2$, $y=u(x-1)$ and
\[
1-x^2-y^2 = (1-x)(1+x-u^2(1-x)).
\]
The strict transform $\wA$ is given by $1+x-u^2(1-x)=0$ and $\omega_0$ is written as
\[
\dfrac{dxdy}{1-x^2-y^2}=\dfrac{-dxdu}{1+x-u^2(1-x)},
\]
hence the assertions are also true on $U_2$. On $U_3$ (resp. $U_4$) the computation is the same as in $U_2$ (resp. $U_1$), with $x$ and $y$ switched throughout.
\end{proof}

The form $\omega_0$ is a regular differential form of top degree on the smooth affine variety $\wX\setminus\wA$. In particular, $\omega_0$ is closed and its restrictions to the irreducible components of the simple normal crossing divisor $\wB\cup E$ vanish for dimension reasons. It follows that $\omega_0$ defines an element of the relative algebraic de Rham cohomology
\[
H^2_\textrm{dR}(\wX\setminus\wA, (\wB\cup E)\setminus \wA) \ = \ M_\textrm{dR}.
\]
(See for instance, Example 3.3.9 in \cite{HM17} or Sections 2.2.7 and 2.2.8 of \cite{BF}.) With abuse of notation, we denote the cohomology class of $\omega_0$ in $M_\mathrm{dR}$ also by $\omega_0$.

Combining with Lemma \ref{Betti} and the integral representation of $\bG$ we obtain the following:

\begin{prop}\label{G is a period of M}
The Catalan constant $\bG$ is a period of the motive $M$. 
\end{prop}

\begin{proof}
Pairing the element $\omega_0$ of the de Rham realization of $M$ with the element $\wD$ of the dual of the Betti realization of $M$ via the period pairing of $M$, in view of $\pi(\wD)=\D$ we get
\[\int_{\widetilde{\Delta}} \omega_0=\int_{\wD}  \dfrac{dxdy}{1-x^2-y^2} \, = \int_{\D} \dfrac{dxdy}{1-x^2-y^2}=\mathbf{G}.\]
(See also Remark \ref{rem: regularity and convergence} below.)
\end{proof}
\begin{rem}\label{rem: regularity and convergence}
Suppose $\omega$ is a meromorphic 2-form on $X$ such that $\pi^\ast\omega$ is regular on $\wX\setminus \wA$ (where $\pi$ is the blow-up map; note that the latter condition implies that $\omega$ is regular on $X\setminus A$). Since $\wX\setminus\wA$ is a 2-dimensional smooth affine variety and $(\wB\cup E)\setminus \wA$ is a simple normal crossing divisor, and the differential form $\pi^\ast\omega$ is a regular top-degree form on $\wX\setminus \wA$ and the topological chain $\wD$ in $\wX^\an\setminus \wA^\an$ is differentiable with boundary in $\wB^\an\cup E^\an$, the period pairing of the class of $\pi^\ast\omega$ in $M_\textrm{dR}$ with $\widetilde{\Delta}\in M^\vee_{\textup{sing}}$ is simply the integral of $\pi^\ast\omega$ over $\wD$ (see \cite[Theorem 5.5.6]{HM17}). By functoriality, this integral is equal to the integral of $\omega$ over $\D$. (In particular, the latter integral converges.)
\end{rem}

\section{The weight filtration of $M$}\label{sec: wt filtration}
The goal of this section is to describe the weight filtration $W_\bullet M$. Since $\wB \cup E$ is a closed subvariety of $\wX$, and $\wA$ is a smooth closed subvariety of $\wX$ of pure codimension $1$, we have a Gysin long exact sequence for $(\wX, \wB \cup E, \wA)$ in $\mathbf{MM}(\Q)$, which reads:
\begin{equation}\label{mainexact}
\begin{tikzcd}[column sep = tiny, row sep = small]
\cdots \longrightarrow H^0(\wA, \wA\cap(\wB\cup E))(-1) &&& \\
\ \ \quad \longrightarrow H^2(\wX, \wB \cup E)  \ar[r] & H^2(\wX \setminus \wA, (\wB \cup E) \setminus \wA) \ar[rr, " \delta "] & & H^1(\wA, \wA \cap (\wB \cup E))(-1) \rightarrow \\ 
 &M\ar[u,equal] & & H^3(\wX, \wB \cup E) \rightarrow \cdots
\end{tikzcd}
\end{equation}
(See \S 4.8 of \cite{FJ} and in particular, Theorem 4.8.3 therein, for the Gysin sequence in the setting of exponential motives, which implies the existence of the sequence in the settings of Nori and absolute Hodge motives.) The map $\delta$ is the connecting morphism. We will be using the sequence above to compute $W_\bullet M$.

\begin{lem}\label{excision lemma}
For any integer $n \geq 0$, the map $\pi \colon \wX \to X$ induces isomorphisms
\begin{align}
&\pi^*\colon H^n(X, B) \xrightarrow{\simeq} H^n(\wX, \wB \cup E)\label{isom2}\\
&\pi^*\colon H^n(A, A \cap B) \xrightarrow{\simeq} H^n(\wA, \wA \cap (\wB \cup E)). \label{isom1}
\end{align} 
\end{lem}

\begin{proof}
The fact that these maps exist in $\mathbf{MM}(\Q)$ is by functoriality.
The map $\pi$ gives an isomorphism between the pairs $(\wA, \wA \cap (\wB \cup E))$ and $(A, A \cap B)$, so the fact that \eqref{isom1} is an isomorphism is clear. We need to prove that the map \eqref{isom2} is an isomorphism. It suffices to show that this map is an isomorphism after passing to the Betti realization, or equivalently, that the map
$$\pi_* \colon H_{n,\textup{sing}}(\wX, \wB \cup E) \to H_{n,\textup{sing}}(X, B)$$
is an isomorphism. There is a commutative diagram
\[
\begin{tikzcd}
   H_{n,\textup{sing}}(\wX, \wB \cup E) \ar[d, "\pi_\ast"] \ar[r, "\simeq"] & H_{n,\textup{sing}}(\wX/(\wB \cup E), \{pt\}) \ar[d, equal] \\
   H_{n,\textup{sing}}(X, B) \ar[r, "\simeq"] & H_{n,\textup{sing}}(X/B, \{pt\}),
\end{tikzcd}
\]
where $pt$ stands for a singleton, and the notation $H_{n,\textup{sing}}(X/B, \{pt\})$ means the relative rational singular homology of the pair $(X^{\an}/B^{\an}, B^{\an}/B^{\an})$, even though $X/B$ is not a variety. The notation $H_{n,\textup{sing}}(\wX/(\wB \cup E), \{pt\})$ is to be interpreted similarly. The horizontal maps are induced by quotient maps and are isomorphisms by a consequence of the classical excision theorem (see e.g., \cite[Proposition 2.22]{Hatcher}), and the vertical maps are induced by $\pi$. The desired conclusion thus follows.
\end{proof}

Throughout the paper, we will use the isomorphisms of Lemma \ref{excision lemma} to identify $H^n(X,B)$ (resp. $H^n(A, A \cap B)$) with $H^n(\wX,\wB\cup E)$ (resp. $H^n(\wA, \wA \cap (\wB \cup E))$).

The topological chain $\D$ in $X^\an$ defines an element of $H_{2,\textup{sing}}(X,B)$, which we also denote by $\D$.

\begin{lem} \label{Q0} \ 
\begin{enumerate}[label=$(\alph*)$]
\item We have an isomorphism
$$H^2(X, B) \simeq \Q(0). $$
Moreover, $H_{2,\textup{sing}}(X, B)$ is spanned by $\Delta$.
\medskip\par 
\item $H_3(X, B)=0$.
\end{enumerate}
\end{lem}

\begin{proof}
\begin{enumerate}[wide, label=($\alph*$)]
\item In view of the long exact sequence for the pair $(X,B)$ for the relative homology in our category of motives and the fact that $H_2(X)$ and $H_1(X)$ vanish because $X^{\an}$ is contractible, it follows that the connecting morphism $H_2(X,B)\rightarrow H_1(B)$, which is the boundary map in Betti realization, is an isomorphism. Thus to establish part (a) it suffices to show that $H_1(B)$ is isomorphic to $\Q(0)$, and the Betti realization $H_{1,\textup{sing}}(B)$ is generated by $\partial \Delta$.

Set $\ell_j^{\,c}:=B\setminus \ell_j$ for $j=1,2,3$ (see Figure \ref{fig1}). It is not difficult to see that the map
\begin{equation}\label{p1eq15}
H_1(\ell_1^{\,c} \, \cap (\ell_2^{\,c}\cup \ell_3^{\,c}))\rightarrow H_1(\ell_1^{\,c}) \, \oplus H_1(\ell_2^{\,c}\cup \ell_3^{\,c})
\end{equation}
induced by the obvious inclusions is an isomorphism: indeed, after passing to the Betti realizations, the domain and codomain of \eqref{p1eq15} are both 4-dimensional spanned by the classes of a small loop around $P$ in $\ell_3^\an$, a small loop around $(0,0)$ in $\ell_2^\an$, and small loops around $Q$ in $\ell_3^\an$ and $\ell_2^\an$; the Betti realization of the map \eqref{p1eq15} sends each generator to its counterpart. Thus the Mayer--Vietoris sequence for the open cover $B=\ell_1^{\,c} \cup (\ell_2^{\,c} \cup \ell_3^{\,c})$ gives an exact sequence
\begin{align} \label{exact}
0 \to H_{1}(B) &\to H_{0}(\ell_1^{\,c} \cap (\ell_2^{\,c} \cup \ell_3^{\,c})) \to H_{0}(\ell_1^{\,c}) \oplus H_{0}(\ell_2^{\,c} \cup \ell_3^{\,c}) \to H_{0}(B) \to 0.
\end{align}
On noting that $(\ell_1^{\,c} \cap (\ell_2^{\,c} \cup \ell_3^{\,c}))^\an$ has two connected components both of which have $\Q$-rational points, and that $(\ell_1^{\,c})^\an$, $B^\an$ and $(\ell_2^{\,c}\cup\ell_3^{\,c})^\an$ are all connected and have $\Q$-rational points, the sequence reads
$$0 \to H_{1}(B) \to \Q(0)^{\oplus 2} \to \Q(0)^{\oplus 2} \to \Q(0) \to 0. $$
It follows that $H_{1}(B)$ is isomorphic to $\Q(0)$. The element $\partial \Delta$ gives a nonzero element of the Betti realization $H_{1,\textup{sing}}(B)$, completing the argument.
\medskip\par 
\item We have
$$H_{3}(X, B) \simeq H_{2}(B)\simeq 0,$$
where the first isomorphism is by the long exact sequence of the relative homology for the pair $(X,B)$. The second isomorphism can be seen from the Mayer--Vietoris long exact sequence for the same open cover of $B$ as in the argument for part $(a)$, in view of injectivity of \eqref{p1eq15} and the fact, easily verified in the Betti realization, that $H_2(\ell_1^{\,c})$ and $H_2(\ell_2^{\,c}\cup \ell_3^{\,c})$ vanish.\qedhere
\end{enumerate}
\end{proof}

By \eqref{mainexact} and Lemmas \ref{excision lemma} and \ref{Q0}, the motive $M$ fits in the following exact sequence: 
\begin{equation*}
\begin{tikzcd}[column sep=small, row sep=small]
 H^0(A, A\cap B)(-1) \ar[d,equal] \ar[r] & H^2(X,B) \ar[r] & M \ar[r, "\delta"] & H^1(A, A \cap B)(-1) \ar[d,equal] \ar[r] & 0\\
H^0(\wA, \wA\cap(\wB\cup E))(-1) &\Q(0) \ar[u, phantom, sloped, "\simeq"] && H^1(\wA, \wA \cap (\wB\cup E))(-1) & 
\end{tikzcd}
\end{equation*}
We now turn our focus to $H^1(A, A \cap B)(-1)$ and $H^0(A, A \cap B)(-1)$. The long exact sequence for the relative homology of the pair $(A,A\cap B)$ gives an exact sequence
\begin{align}\label{p1eq4}
0=H_1(A \cap B) \to H_1(A) \to H_1(A, A \cap B) \to H_0(A \cap B) \xrightarrow{(*)} H_0(A) \to H_0(A, A \cap B) \to 0.
\end{align}
Note that $A^\an$ is connected and $A \cap B$ consists of $4$ points defined over $\Q$, so that 
\[
H_0(A\cap B) \simeq \Q(0)^{\oplus 4} \hspace{.3in} \text{and} \hspace{.3in} H_0(A) \simeq \Q(0),
\]
with Betti realizations given by 
\begin{align*}
H_{0, \textup{sing}}(A \cap B)= \langle [P], [Q], [-P], [-Q]\rangle_\Q \hspace{.3in} \text{and} \hspace{.3in} H_{0, \textup{sing}}(A)=\langle[\,pt\,]\rangle_\Q,
\end{align*}
where $[\,pt\,]$ refers to the class of a  singleton; here and elsewhere in the paper, we use the notation $\langle \ \ \rangle_\Q$ for the $\Q$-span. The map
$$(*)_{\textup{sing}}  \colon H_{0, \textup{sing}}(A \cap B) \to H_{0, \textup{sing}}(A)$$
sends each of $[\pm P]$, $[\pm Q]$ to $[\,pt\,]$}. It follows that $(*)$ is surjective with kernel $\Q(0)^{\oplus 3}$, generated in the Betti realization by $[Q]-[P]$, $[-P]-[P]$ and $[-Q]-[P]$. 
Thus $H^0(A,A\cap B)$ is zero, and moreover, we have a short exact sequence
$$0 \to H_1(A) \to H_1(A, A \cap B) \to \Q(0)^{\oplus 3} \to 0. $$
Taking duals, we get a short exact sequence 
\begin{align} \label{mainexact3}
0 \to \Q(0)^{\oplus 3} \to H^1(A, A \cap B) \to H^1(A) \to 0.
\end{align}
For future referencing, we record our conclusions:
\begin{lem}\label{lem: H^1(A,A cap B)} \ 
\begin{enumerate}[label=$(\alph*)$]
\item The sequence \eqref{mainexact} gives a short exact sequence
\begin{equation}\label{mainexact2}
\begin{tikzcd}[column sep=small, row sep=small]
0 \ar[r] & H^2(X,B) \ar[r] & M \ar[r, "\delta"] & H^1(A, A \cap B)(-1) \ar[r] & 0.\\
&\Q(0) \ar[u, phantom, sloped, "\simeq"] && & 
\end{tikzcd}
\end{equation}
\item $H^1(A, A \cap B)(-1)$ is an extension of $H^1(A)(-1)$ by $\Q(-1)^{\oplus 3}$ in $\mathbf{MM}(\Q)$.
\end{enumerate}
\end{lem}

Let $K=\Q(i)\subset\C$. After base change to $K$, we have an isomorphism  
\begin{equation}\label{p1eq6}
A_K=\operatorname{Spec}{K}[x, y]/(1-x^2-y^2) \xrightarrow{\simeq} \mathbb{G}_{m,K} = \operatorname{Spec} K[z, z^{-1}]; \quad (x, y) \mapsto z=x+yi
\end{equation}
over $K$, so that
\[H^1(A_K) \simeq H^1(\mathbb{G}_{m,K}) \simeq \Q(-1)\]
in $\mathbf{MM}(K)$. Since the weight filtration behaves well under base change from $\Q$ to $K$ (as the weight filtration is determined by the Hodge realization), it follows that the motive $H^1(A)$ (in $\mathbf{MM}(\Q)$) is pure of weight 2. In view of Lemma \ref{lem: H^1(A,A cap B)} and the fact that morphisms of motives are strictly compatible with the weight filtration, we obtain that
\begin{align}\label{p1eq8}
\begin{split}
    W_0M &= H^2(\wX,\wB\cup E)=H^2(X,B)\simeq\Q(0)\\
    W_2M &= \delta^{-1}\left(W_2\left[H^1(A,A\cap B)(-1)\right]\right) \\
    &= \delta^{-1}\biggm(\ker\left[H^1(A,A\cap B)\rightarrow H^1(A)\right]\otimes\ \Q(-1) \, \simeq \, \Q(-1)^{\oplus 3}\bigg)\\
    W_4M &=M 
\end{split}
\end{align}
and
$$\operatorname{Gr}^W(M) \simeq  \Q(0) \oplus \Q(-1)^{\oplus 3} \oplus H^1(A)(-1)$$
(with weights 0,2,4). In particular, $M$ is $5$-dimensional.

\begin{rmk}
The motive $H^1(A)$ is not isomorphic to $\Q(-1)$ in $\mathbf{MM}(\Q)$, as the periods of $H^1(A)$ are in $\pi \Q$ rather than in $\pi i \Q$ (see the period computations of the next section). Thus $M$ is not a mixed Tate motive over $\Q$. (By a mixed Tate motive over $\Q$ here, we simply mean an object of $\mathbf{MM}(\Q)$ that is an iterated extension of the $\Q(n)$.)
\end{rmk}

\section{The de Rham and Betti realizations}\label{sec: Betti-dR realization}
\numberwithin{thm}{subsection}
In this section, we explicitly describe the de Rham realization of the motive $M$ defined in \S \ref{construction} and the Betti realization of its dual $M^\vee$, starting with the latter. We then use these to compute the periods of $M$.

\subsection{The Betti realization}\label{sebsection: Betti realization}
Dualizing \eqref{mainexact2} and taking the Betti realization, we obtain a short exact sequence
\begin{equation}\label{p1eq3}
\begin{tikzcd}
0\ar[r] &   H_{1,\textup{sing}}(A, A \cap B)(1) \ar[r, "\delta^\vee_{\textup{sing}}"] & M^\vee_{\textup{sing}} \ar[r] & H_{2,\textup{sing}}(X, B) \ar[r] & 0.
\end{tikzcd}
\end{equation}

First, we determine a basis of $H_{1,\textup{sing}}(A, A \cap B)$. By \eqref{p1eq4},
we have a short exact sequence 
$$0 \to H_{1, \textup{sing}}(A) \to H_{1, \textup{sing}}(A, A \cap B) \xrightarrow{\partial} \ker \biggm(H_{0, \textup{sing}} (A \cap B) \xrightarrow{(\ast)_{\textup{sing}}} H_{0, \textup{sing}}(A) \biggm) \to 0, $$
where the connecting homomorphism $\partial$ is given by the boundary map. Since $A^\an\simeq \C^\times$ via $(x,y)\mapsto x+iy$, we have $H_{1, \textup{sing}}(A) \simeq \Q \cdot \gamma$, where $\gamma$ is the unit circle in $\R^2$ around $(x, y)=(0, 0)$ oriented positively. For a point $R \in \{- P, Q,-Q\}$, let $\alpha_{_{P,R}}$ be the counter-clockwise arc from $P$ to $R$ on the unit circle $\gamma$. On recalling the description of the kernel of $(\ast)_{\textup{sing}}$ from \S \ref{sec: wt filtration} (see the paragraph before Lemma \ref{lem: H^1(A,A cap B)}), $\gamma$ together with the $\alpha_{_{P,R}}$ with $R\in \{- P, Q,-Q\}$ form a basis of $H_{1, \textup{sing}}(A, A \cap B)$. It will be slightly more convenient to adjust this basis to the basis $\{\gamma,\alpha_1,\alpha_2,\alpha_3\}$, where
\begin{align}\label{p1eq9}
&\alpha_1:=-\alpha_{_{P,Q}}+\alpha_{_{P,-Q}}-\gamma \, , \quad  \alpha_2:=-\alpha_{_{P,-P}}+\frac{1}{2} \, \gamma \, , \quad \alpha_3:=-2\alpha_{_{P,Q}}+2\alpha_{_{P,-P}}-2\alpha_{_{P,-Q}}+\gamma.
\end{align}

On the other hand, the map $M^\vee_{\textup{sing}}\rightarrow  H_{2,\textup{sing}}(X, B)$ in \eqref{p1eq3} is induced by the blow-up map. By Lemma \ref{Q0}, $H_{2,\textup{sing}}(X, B)$ is 1-dimensional generated by the simplex $\Delta$, which has a lift $\wD$ in $M^{\vee}_{\textup{sing}}$ (see Lemma \ref{Betti} and the comments before and after it). In view of the sequence \eqref{p1eq3}, we obtain the statement in part $(a)$ below.

\begin{lem}\label{lem: Betti realization} \ 
\begin{enumerate}[label=$(\alph*)$]
\item A basis of $M^\vee_{\textup{sing}}$ is formed by
\begin{align} \label{baseB}
\wD, \  \gamma \otimes 2 \pi i, \ \alpha_j \otimes 2 \pi i \ \ \ \text{for $j = 1,2,3$},
\end{align}
where we have considered $\gamma \otimes 2 \pi i$ and $\alpha_j \otimes 2 \pi i$ ($j=1,2,3$) as elements of $M^\vee_{\textup{sing}}$ via the injection $\delta^\vee_{\textup{sing}}$ in \eqref{p1eq3}. 
\item The weight filtration on $M^\vee_{\textup{sing}}$ is given in terms of the basis \eqref{baseB} by
\begin{align*}
W_{-4}(M^\vee_{\textup{sing}}) &= \langle \gamma \otimes 2 \pi i\rangle_\Q\\
W_{-2}(M^\vee_{\textup{sing}}) &= \langle \gamma \otimes 2\pi i,\  \alpha_1 \otimes 2 \pi i, \alpha_2 \otimes 2 \pi i, \alpha_3 \otimes 2 \pi i \rangle_\Q\\
W_{0}(M^\vee_{\textup{sing}}) &= M^\vee_{\textup{sing}}.
\end{align*}
\end{enumerate}
\end{lem}
\medskip\par 
\begin{proof}
We have already established part $(a)$. Dualizing \eqref{mainexact2}, the weight filtration on $M^\vee$ is given by $W_{-4}(M^\vee)=W_{-4}[H_1(A,A\cap B)(1)]=H_1(A)(1)$, $W_{-2}(M^\vee)=H_1(A,A\cap B)(1)$, and $W_0(M^\vee)=M^\vee$. The assertion in part $(b)$ follows. 
\end{proof}

\subsection{The de Rham realization}\label{subsection: dR realization}
We now compute an explicit basis of $M_\textrm{dR}$. Taking the de Rham realization of \eqref{mainexact2}, we have a short exact sequence
\begin{equation} \label{p1eq5}
\begin{tikzcd}[row sep=small]
0 \ar[r] & H_\textrm{dR}^2(X,B) \ar[r] & M_\textrm{dR}  \ar[r,"\operatorname{Res}_{\wA}"]  &H_\textrm{dR}^1(A, A \cap B)(-1) \ar[d,equal] \ar[r] & 0,\\
&\Q(0)_\textrm{dR}=\Q \ar[u, phantom, sloped, "\simeq"] & &H_\textrm{dR}^1(A, A \cap B) & 
\end{tikzcd}
\end{equation}
where $\operatorname{Res}_{\wA}$ is the Poincar\'{e} residue map along $\wA$. (Recall that we have identified $H^1(\wA,\wA\cap(\wB\cup E))$ and $H^1(A,A\cap B)$ via the isomorphism \eqref{isom1}.)
\medskip\par 

We start by finding a basis of $H^1_{\rm dR}(A, A \cap B)$.
\begin{notation} \label{baromega} 
Let $\overline{\omega}_j$ ($0\leq j\leq 3$) be the following differential forms on $A$ and $\wA$:
\begin{alignat*}{2}
&\overline{\omega}_0=\dfrac{dx}{2y}=\dfrac{-dy}{2x}, \quad \quad &&\overline{\omega}_1=x \,\overline{\omega}_0=-\dfrac12 dy, \\
&\overline{\omega}_2=y \,\overline{\omega}_0=\dfrac12 dx, \quad \quad &&\overline{\omega}_3=xy \, \overline{\omega}_0=\dfrac12 xdx.  
\end{alignat*}
\end{notation}
Note that since $x^2+y^2=1$ on $A$, we may write $\overline{\omega}_0$ also as $\frac{1}{2}(y\,dx-x\,dy)$. For each $j$, the differential form $\overline{\omega}_j$ is a top-degree regular form on $A$, so that it defines an element of $H^1_{\rm dR}(A, A \cap B)$. With abuse of notation, we denote the elements of $H^1_{\rm dR}(A, A \cap B)$ and $H^1_{\rm dR}(A)$ defined by $\overline{\omega}_j$ also by $\overline{\omega}_j$.

Under the isomorphism \eqref{p1eq6}, the differential forms $\ow_j$ are given by
\begin{align}\label{p1eq7}
\overline{\omega}_0=\dfrac{idz}{2z}\, , \quad \overline{\omega}_1=\dfrac{i}4 (1+z^{-2})dz\, , \quad \overline{\omega}_2=\dfrac14(1-z^{-2}) dz\, , \quad  \overline{\omega}_3=\dfrac18 z(1-z^{-4}) dz
\end{align}
(as $z+\dfrac{1}{z}=2x$ and $z-\dfrac{1}{z}=2iy$). The period matrix of $\ow_j$ ($j=0,\ldots,3$) with respect to the basis $\alpha_1,\alpha_2,\alpha_3,\gamma$ of $H_{1,\textup{sing}}(A,A\cap B)$ (see \eqref{p1eq9}) is computed to be\footnote{The reason for the change of basis \eqref{p1eq9} was to obtain a diagonal period matrix here.}
\begin{equation*}
 \begin{blockarray}{rcccccccc}
 & & \ow_1 & \ow_2& \ow_3 & \ow_0   \\
 \begin{block}{rr(ccccccc)}
\alpha_1 && 1 & 0& 0& 0  \\
\alpha_2&& 0 & 1& 0& 0 \\
\alpha_3  && 0&0 & 1 & 0      \\
\gamma & & 0& 0& 0&  -\pi   \\ 
 \end{block}
 \end{blockarray} \ . 
 \end{equation*}
The entries are calculated using computations such as
\[\int_{\alpha_{_{P,Q}}} \ow_0= \int_1^i \dfrac{i \,dz}{2z}= \dfrac{i}{2} \cdot \dfrac{\pi i}{2}=-\frac{\pi}4.\]
Since the period matrix above is invertible, it follows that the $\ow_j$ ($j=0,1,2,3$) form a basis of $H^1_\textrm{dR}(A,A\cap B)$. We record the conclusion:

\begin{lem}\label{period matrix of H^1(A,A cap B)}
The elements $\ow_1,\ow_2,\ow_3,\ow_0$ form a basis of $H^1_{\rm dR}(A,A\cap B)$. The period matrix of $H^1(A,A\cap B)$ with respect to the ordered basis $(\ow_1,\ow_2,\ow_3,\ow_0)$ of $H^1_{\rm dR}(A,A\cap B)$ and the ordered basis $(\alpha_1, \alpha_2, \alpha_3, \gamma)$ of $H_{1,\textup{sing}}(A,A\cap B)$ (see \eqref{p1eq9}) is the diagonal matrix with diagonal entries $1,1,1,-\pi$, in that order.
\end{lem}

\begin{rmk}\label{rem: wt filt of H^1_dR(A, A cap B)}
The de Rham realization of \eqref{mainexact3} is given by 
$$0 \to \Q(0)^{\oplus 3}_{\rm dR} \to H^1_{\rm dR}(A, A \cap B) \xrightarrow{(\dagger)} H^1_{\rm dR}(A) \to 0. $$
In view of \eqref{p1eq7} it is clear that $\ow_1, \ow_2$ and $\ow_3$ generate the kernel of $(\dagger)$ and $\ow_0$ maps to a generator of $H^1_{\rm dR}(A)$ (also denoted by $\ow_0$). Thus 
\[
\bigm(W_0H^1(A,A\cap B)\bigm)_\textrm{dR} \ = \ \langle\ow_1,\ow_2,\ow_3\rangle_\Q.
\]
\end{rmk}

\begin{rmk} \label{splitK}
The period matrix of $H^1(A,A\cap B)$ suggests in view of Grothendieck's period conjecture that the extension $H^1(A,A\cap B)$ of $H^1(A)$ by $\Q(0)^{\oplus 3}$ should be split in $\mathbf{MM}(\Q)$. Note that the fact that $H^1(A,A\cap B)$ splits in $\mathbf{MM}(K)$ (with $K=\Q(i)$ as before) is clear to the expert reader, at least assuming that $\operatorname{Ext}^1_{\mathbf{MM}(K)}(\Q(0),\Q(1))$ is given by Kummer motives (as it is expected), since $A\cap B$ translates under \eqref{p1eq6} to the torsion subgroup $\mu_4$ of $\mathbb{G}_{m,K}$ (recall that throughout, we are working with motives with rational coefficients). We will see shortly that $H^1(A,A\cap B)$ is indeed also split in $\mathbf{MM}(\Q)$. See Lemma \ref{splitQ}.
\end{rmk}

We now lift the elements $\ow_0, \ldots, \ow_3$ of $H_{\rm dR}^1(A, A \cap B)$ to elements of $M_{\rm dR}$, as follows. 
\begin{notation}\label{notation: def of omega_j}
From now on, for simplicity we set 
\[g:=1-x^2-y^2.\]
Let $\omega_j$ ($1\leq j\leq 3$) be the following differential forms on $X$ and $\wX$ (recall the concrete description of $\wX$ in \S \ref{construction}):
\[ \omega_1:=x \omega_0=\dfrac{xdxdy}{g}, \quad \omega_2:=y\omega_0=\dfrac{ydxdy}{g}, \quad \omega_3:=xy \omega_0=\dfrac{xydxdy}{g}, \]
where
\[
\omega_0=\dfrac{dxdy}{g}
\]
as in \eqref{p1eq2}.
\end{notation}
By Lemma \ref{deRham}, the $\omega_j$ are regular on $\wX \setminus \wA$ with logarithmic singularities along $\wA$. In particular, the Poincar\'{e} residue of each $\omega_j$ along $\wA$ is well-defined.

\begin{lem}\label{lem: Res of omega_j}
For $j=0, \ldots, 3$, we have 
$$\operatorname{Res}_{\wA} (\omega_j)=\ow_j,  $$
where $\ow_j$ is the $1$-form defined in Notation  \ref{baromega} and $\operatorname{Res}_{\wA}$ is the Poincar\'e residue map along $\wA$. (The equality above is an equality of differential forms on $\wA$, or on $A$ via the isomorphism $\wA\cong A$ given by the blow-up map.)
\end{lem}

\begin{proof}
We verify the equality outside $xy=0$. The strict transform $\wA$ is then defined by $g=0$. Since $dg=-2x dx -2y dy$, we have
\begin{align*}
\omega_0=\dfrac{dxdy}{g}= \dfrac{dg}{g} \wedge \dfrac{dx}{2y} \, .
\end{align*}
Thus $\operatorname{Res}_{\wA} (\omega_0)$ and $\ow_0$ coincide on $\wA\setminus \{xy=0\}$. Since both $\operatorname{Res}_{\wA} (\omega_0)$ and $\ow_0$ are regular on $\wA$, it follows that $\operatorname{Res}_{\wA} (\omega_0)=\ow_0$ on $\wA$.
The assertions for $\omega_1, \omega_2$ and $\omega_3$ are  checked by inserting $x$, $y$ and $xy$ in appropriate places of the computation.
\end{proof}

Being a regular differential form of $\wX\setminus\wA$ of top degree, each $\omega_j$ defines an element of $M_{\rm dR}$. With abuse of notation, we denote this element also by $\omega_j$. The equation of Lemma \ref{lem: Res of omega_j} descends to the same equation with $\omega_j$ and $\ow_j$ respectively interpreted as elements of $M_\textrm{dR}$ and $H^1_\textrm{dR}(A,A\cap B)$. Combining with Lemma \ref{period matrix of H^1(A,A cap B)} and in view of the sequence \eqref{p1eq5}, it follows that to complete $\{\omega_0,\omega_1, \omega_2,\omega_3\}$ to a basis of $M_{\mathrm{dR}}$ it suffices to add a generator of the line $H^2_{\rm dR}(X, B) \subset M_{\rm dR}$. 

Set 
\begin{equation}\label{p1eq13}
    \eta := 2dxdy,
\end{equation}
which (with abuse of notation) is a regular differential form on $\wX$ as well as on $X$. It defines an element of $M_\textrm{dR}$, also denoted by $\eta$, whose Poincar\'{e} residue along $\wA$ vanishes. Hence $\eta$ is in $H^2_{\rm dR}(X, B) \subset M_{\rm dR}$. Finally, its period pairing with the element $\wD$ of $M^\vee_{\textup{sing}}$ is given by
$$\int_{\wD} \eta = \int_{\Delta} \eta=1. $$
It follows that $\eta$ is nonzero in $M_\textrm{dR}$ and hence it generates $H^2_{\rm dR}(X, B)$. 
\medskip\par
Combining with \eqref{p1eq8} and Remark \ref{rem: wt filt of H^1_dR(A, A cap B)}, we have proved the following result:
\begin{lem}\label{lem: dR realization}
A basis of $M_{\rm dR}$ is given by 
\begin{align} \label{basedR}
\{\eta, \omega_0, \ldots, \omega_3\} 
\end{align}
(where the $\omega_j$ and $\eta$ are defined in Notation \ref{notation: def of omega_j} and equation \eqref{p1eq13}). Moreover, 
$$W_0 M_{\rm dR} = \langle \eta \rangle_\Q, \quad W_2 M_{\rm dR}=\langle \eta, \omega_1, \omega_2, \omega_3 \rangle_\Q, \quad W_4M_{\rm dR}=M_{\rm dR}.$$
\end{lem}

\subsection{The period matrix of $M$}
In the previous two subsections we found bases of $M^\vee_{\textup{sing}}$ and $M_{\rm dR}$. We now compute the period matrix of $M$ with respect to these bases.

\begin{prop}\label{prop: period matrix}
The period matrix of $M$ with respect to the basis of $M^\vee_{\textup{sing}}$ given in Lemma \ref{lem: Betti realization} and the basis of $M_\textup{dR}$ given in Lemma \ref{lem: dR realization} is 
\begin{equation}\label{eq per matrix  of M}
 \begin{blockarray}{rcccccccc}
 & &  \eta & \omega_1 & \omega_2& \omega_3 & \omega_0   \\
 \begin{block}{rr(ccccccc)}
\wD \quad \quad   & & 1& (\log 2)/2  &  (\log 2)/2 &  (1- \log 2)/4 & \bG \\
\alpha_1\otimes 2\pi i & &  0 & 2 \pi i & 0& 0& 0  \\
\alpha_2\otimes 2\pi i & & 0 & 0& 2 \pi i & 0& 0  \\
 \alpha_3\otimes 2\pi i && 0 & 0&0 & 2 \pi i & 0  \\
 \gamma\otimes 2\pi i & & 0& 0 & 0& 0&  -\pi \cdot (2 \pi i)  \\ 
 \end{block}
 \end{blockarray} \ .
 \end{equation}
\end{prop}

\begin{proof}
That the (1,5) entry is $\bG$ is by the integral representation of $\bG$ (see Proposition \ref{G is a period of M}). The bottom right $4\times 4$ submatrix comes from the computation of periods of $H^1(A,A\cap B)$ (Lemma \ref{period matrix of H^1(A,A cap B)}). That the (1,1) entry is 1 is clear. The fact that the remaining entries of the first column are zero can be seen, for instance, by weight considerations: the elements $\alpha_j\otimes 2\pi i$ and $\gamma\otimes 2\pi i$ are in $W_{-1}(M^\vee_{\textup{sing}})$, and hence they send $W_0M_\textrm{dR}$ to zero. Thus the only new claims are about how $\omega_j$ for $j=1,2,3$ pairs with $\wD$. These are done by direct computations. For example,
\begin{align*}
\int_{\wD} \omega_1=\int_{\Delta} \dfrac{x\,dxdy}{1-x^2-y^2}=\int_0^1\left(\int_0^{1-y} \dfrac{x\,dx}{1-x^2-y^2}\right) dy =\frac{\log 2}{2}
  \end{align*}
(see also Remark \ref{rem: regularity and convergence}). The other integrals are similarly calculated.  
\end{proof}

\begin{rmk}\label{rem: why log2}
The top left $4\times 4$ submatrix of \eqref{eq per matrix  of M} is a period matrix of $W_2M$, which is an extension of $\Q(-1)^{\oplus 3}$ by $\Q(0)$. The existence of $\log 2$ among the periods and its irrationality show that this extension is nonsplit in $\mathbf{MM}(\Q)$, as its Betti--de Rham realization is nonsplit. The existence of $\log 2$ among the periods suggests that $W_2M$ (and hence $M$) is ramified at 2.
\end{rmk}

\section{The motive of the Catalan constant} \label{Catalanmotive}
\numberwithin{thm}{section}
In this section, we define (what we will call) the motive of the Catalan constant as a particular 2-dimensional submotive of $M$. Recall that we have a diagram
\begin{equation}\label{p1eq10}
\begin{tikzcd}[column sep = small]
&&& 0 \ar[d] & \\
&&&\Q(-1)^{\oplus 3} \ar[d] & \\
0 \ar[r] & H^2(X, B)    \ar[r] & M  \ar[equal,d]  \ar[r, "\delta"] & H^1(A, A \cap B)(-1) \ar[d] \ar[r] &  0  \\
& \Q(0) \ar[u, phantom, sloped, "\simeq"]& H^2(\wX \setminus \wA, (\wB \cup E) \setminus \wA)  & H^1(A)(-1) \ar[d] & \\
&&&0  &
\end{tikzcd}
\end{equation}
in $\mathbf{MM}(\Q)$, with exact horizontal and vertical sequences.

The pair $(A, A \cap B)=(A, \{\pm P,\pm Q\})$ has a natural action of the dihedral group $D_4$. In particular, the automorphism $\sigma$ of $X=\operatorname{Spec}\Q[x,y]$ defined by 
$$\sigma \colon (x, y) \mapsto (-y, x) $$
gives an automorphism of order 4 of the pair $(A, A \cap B)$, and hence an automorphism of the motive $H^1(A, A \cap B)$, which we also denote by $\sigma$. There is thus a motive 
$$H^1(A, A \cap B)^{\sigma} := \ker\bigm(\operatorname{Id} -\sigma \in \operatorname{End}_{\mathbf{MM}(\Q)}(H^1(A, A \cap B))\bigm)$$
in $\mathbf{MM}(\Q)$. The action of $\sigma$ on the de Rham realization of $H^1(A, A \cap B)$ is given in terms of the explicit basis $\{\ow_j: 0\leq j\leq 3\}$ of $H^1_\textrm{dR}(A, A \cap B)$ (see Notation \ref{baromega}) by
$$\sigma(\ow_0)=\ow_0, \quad \sigma(\ow_1)=-\ow_2, \quad \sigma(\ow_2)=\ow_1, \quad \sigma(\ow_3)=-\ow_3. $$
It follows that $H^1(A, A \cap B)^{\sigma}$ is $1$-dimensional with de Rham realization
\begin{equation}\label{p1eq14}
{H^1_\textrm{dR}(A, A \cap B)^{\sigma}} \, = \, \langle \ow_0 \rangle_\Q \, .
\end{equation}

\begin{lemma} \label{splitQ}
The sequence
$$0 \to \Q(0)^{\oplus 3} \to H^1(A, A \cap B) \to H^1(A) \to 0$$
(or equivalently, the vertical sequence in \eqref{p1eq10}) splits in $\mathbf{MM}(\Q)$. 
\end{lemma}
\begin{proof}
Consider the diagram 
\[
\begin{tikzcd}
&&  H^1(A, A \cap B)^{\sigma}  \ar[d, hookrightarrow]  \ar[dr, "(\ddagger)"] & \\
0 \ar[r] & \Q(0)^{\oplus 3} \ar[r] & H^1(A, A \cap B) \ar[r]& H^1(A) \ar[r] & 0.\\
\end{tikzcd}
\]
The map $(\ddagger)$ (defined by the commutativity of the diagram) is an isomorphism since it is an isomorphism after passing to the de Rham realizations (it sends $\ow_0$ in $H^1_\textrm{dR}(A,A\cap B)^\sigma$ to its namesake in $H^1_\textrm{dR}(A)$). This gives the desired conclusion. 
\end{proof}

\begin{dfn}
The \textit{motive of the Catalan constant}, denoted by $C$, is the pullback of $M$ under the inclusion $H^1(A, A \cap B)^{\sigma} (-1) \subset H^1(A, A \cap B)(-1)$. That is, $C$ is the subobject of $M$ fitting in the commutative diagram
\[
\begin{tikzcd}
0 \ar[r]& H^2(X, B) \ar[r] \ar[equal,d] & C \ar[r] \ar[d, hookrightarrow] & H^1(A, A \cap B)^{\sigma}(-1) \ar[r] \ar[d, hookrightarrow] & 0 \\
0 \ar[r] & H^2(X, B) \ar[r] 
\ar[r] & M \ar[r]& H^1(A, A \cap B)(-1) \ar[r]& 0.
\end{tikzcd}
\]
\end{dfn}
\medskip\par 
Note that $C$ is a 2-dimensional object of $\mathbf{MM}(\Q)$ with 
\[
\operatorname{Gr}^W(C)  \ = \ H^2(X,B) \, \oplus H^1(A, A \cap B)^{\sigma}(-1) \ \simeq \ \Q(0) \oplus H^1(A)(-1).
\]
The motive $C$ is not a mixed Tate motive over $\Q$, however, its extension of scalars to $K=\Q(i)$ is mixed Tate with
\[
\operatorname{Gr}^W(C_K) \ \simeq \ \Q(0)  \oplus \Q(-2) \quad \text{in} \quad \mathbf{MM}(K).
\]
\medskip\par 
In view of \eqref{p1eq14} and \S \ref{subsection: dR realization} (particularly, Lemmas \ref{lem: dR realization} and \ref{lem: Res of omega_j}), the de Rham realization of $C$ has a basis given by $\eta$ and $\omega_0$. There is a natural surjection
\[
M^\vee_{\textup{sing}}\rightarrow C^\vee_{\textup{sing}}.
\]
With abuse of notation, we denote the images of $\gamma\otimes 2\pi i$ and $\wD$ under this map by the same symbols. By functoriality of the period pairing, the elements $\eta,\omega_0$ of $C_\textrm{dR}$ pair with the elements $\gamma\otimes 2\pi i$ and $\wD$ of $C^\vee_{\textup{sing}}$ to give the period matrix
\begin{equation}\label{p2eq12}
 \begin{blockarray}{rcccccccc}
 & & \eta & \omega_0  \\
 \begin{block}{rr(ccccccc)}
\wD \quad \quad & & 1 & \mathbf{G}  \\
\gamma \otimes 2 \pi i & & 0& -\pi \cdot (2 \pi i)   \\  
 \end{block} 
 \end{blockarray} \ . 
 \end{equation}
It follows that $\gamma\otimes 2\pi i$ and $\wD$ give a basis of $C^\vee_{\textup{sing}}$.

\begin{rmk}
We do not expect $\sigma$ to extend to an automorphism of the motive $M$, as even though $\sigma$ acts on $X\setminus A$, it does not keep $B$ invariant.
\end{rmk}

\begin{rmk}\label{rem: automorphism tau}
The only nontrivial element of $D_4$ that acts on the pair $(X\setminus A, B\setminus A)$ is the reflection
\[
\tau \colon (x,y)\mapsto (y,x).
\]
(Note that $\sigma$ and $\tau$ are the `standard' generators of $D_4$.) The automorphism $\tau$ extends to an automorphism of the pair $(\wX\setminus \wA, (\wB\cup E)\setminus \wA)$ given by
\[
\left((x, y), [\alpha_1 : \beta_1], [\alpha_2 : \beta_2]\right) \mapsto \left((y, x), [\beta_2 : \alpha_2], [\beta_1 : \alpha_1] \right),
\]
and hence gives an automorphism of the motive $M$, which we also denote by $\tau$. We cannot fully extract the motive $C$ from $M$ using $\tau$. Indeed, the motive $M$ has an eigenspace decomposition
$$M=M^\tau  \oplus M^{\tau=-1},$$
where $M^{\tau =-1}$ (resp. $M^\tau$) denotes the eigenspace on which $\tau$ acts as $-1$ (resp. as identity). Direct computations in the de Rham realization give
$$\tau(\eta)=-\eta, \quad \tau(\omega_0)=-\omega_0, \quad \tau(\omega_1)=-\omega_2, \quad \tau(\omega_2)=-\omega_1, \quad \tau(\omega_3)=-\omega_3, $$
so that
$$M^{\tau=-1}_{\rm dR}=\langle \eta, \omega_0, \omega_3, \omega_1+\omega_2\rangle_\Q, \quad M^{\tau}_{\rm dR}= \langle \omega_1 -\omega_2 \rangle_\Q. $$
Thus, $\operatorname{dim}M^{\tau=-1}=4$ and $\operatorname{dim}M^{\tau}=1$. On recalling the description of $C_\textrm{dR}$, we have $C \subset M^{\tau=-1}$.
\end{rmk}

\section{Linear forms in 1 and $\bG$ -- I: A motivic proof of Theorem \ref{thm: linear forms}}\label{sec: lin forms I}
\numberwithin{thm}{subsection}
\subsection{Overview}
The main objective of this section is to give a motivic proof of Theorem \ref{thm: linear forms} using our earlier considerations. The proof we give here will be non-constructive, in that it will not give us the coefficients of 1 and $\bG$ in the linear forms.
\medskip\par 
We recall the basic idea here, which follows Brown's ideas \cite{Br16} (further investigated by Dupont \cite{Dupont} and Dupont--Fres\'{a}n \cite{DF}; a very nice exposition of these ideas can also be found in Burgos Gil--Fres\'{a}n's \cite{BF}). Recall that $C$ is a 2-dimensional motive over $\Q$ with
\[
C_\textrm{dR} = \langle \eta,\omega_0\rangle_\Q \ \ \ \text{and} \ \ \ C^\vee_{\textup{sing}} = \langle \widetilde{\Delta}, \gamma\otimes 2\pi i\rangle_\Q
\]
and period matrix as in \eqref{p2eq12}. If $\omega$ is an element of $C_\textrm{dR}$, then
\[
\omega = a_{\omega} \eta + b_{\omega}\omega_0
\]
in $C_\textrm{dR}$ for some $a_{\omega}$ and $b_{\omega}$ in $\Q$. Pairing with $\widetilde{\Delta}$ via the period pairing, we obtain
\[
\int_{\tilde{\Delta}} \omega = a_\omega+ b_\omega \bG
\]
(where with abuse of notation, $\omega$ on the left hand side of the last equation denotes a regular 2-form on $\wX\setminus \wA$ that represents its namesake in $C_{\mathrm{dR}}$). The ultimate hope is that this approach may lead to linear forms in $1$ and $\bG$ of sufficiently good quality that might be useful in proving irrationality of $\bG$. More explicitly, the ultimate hope would be to try to choose a sequence of cohomology classes $\{\omega_{n}\}_{n\in \Z}$ such that
\[
\int_{\tilde{\Delta}} \omega_n \neq 0
\]
and
\[
\left(\int_{\tilde{\Delta}} \omega_n\right)\cdot (\text{the lowest common denominator of $a_{\omega_n}$ and $b_{\omega_n}$}) \ \ \longrightarrow \ 0
\]
as $n\rightarrow \infty$. By a standard elementary argument this would imply irrationality of $\bG$. Of course, to carry out this approach, one likely needs to have an explicit way of computing the coefficients of 1 and $\bG$ in the linear forms. We shall address this problem in \S \ref{sec: lin forms II} and \S \ref{sec: lin forms III}.
\medskip\par 
In \S \ref{sec:integrability} and \S \ref{sec: sym diff forms} below we discuss how one can naturally produce elements of $C_{\dR}$ represented by differential forms on $X\setminus A$. Throughout the rest of the paper, given any $F\in \Q[x,y]$ and integer $t\geq 0$, we set
\begin{equation}\label{p2eq11}
\omega_{F,t}:= \frac{Fdxdy}{g^{t+1}}
\end{equation}
(with $g=1-x^2-y^2$, as before), which can be considered as a meromorphic 2-form on $\widetilde{X}$ or $X$ (the form on $\widetilde{X}$ is the pullback of its counterpart on $X$ under the blow-up map). In \S \ref{sec:integrability} below we study when the meromorphic differential form $\omega_{F,t}$ on $\wX$ is regular on $\wX\setminus \wA$. We will see that this is guaranteed by a weaker version of condition (ii) of Theorem \ref{thm: linear forms}. Once this regularity condition is imposed, $\omega_{F,t}$ represents a cohomology class in $M_\textrm{dR}$ and by our earlier calculation of the periods of $M$, it follows that
\[
\int\limits_\Delta \omega_{F,t} \ \in \ \langle 1, \log 2, \bG\rangle_{\Q}.
\]
In \S \ref{sec: sym diff forms} (Proposition \ref{prop: proof of lin form thm, t=0 case}), we will show that the regularity condition on $\omega_{F,t}$ together with its $\sigma$-invariance will guarantee that $\omega_{F,t}$ represents an element of $C_\dR$ and hence the integral above will belong to $\langle 1,\bG\rangle_\Q$. This will establish a more natural and (only seemingly, see Remark \ref{rem: ver 1 implies ver 2}) more general version of Theorem \ref{thm: linear forms}, in which condition (ii) is replaced by a regularity condition (see Theorem \ref{thm: linear forms, ver2}). The final subsection will start our preparations for the last objective of the paper, which will be to explicitly compute the coefficients of 1 and $\bG$ in our linear forms.

\subsection{Regularity}\label{sec:integrability}
Let $F\in \Q[x,y]$ be a nonzero polynomial. Then $F$ can be considered as a regular function on $\widetilde{X}$. Recall that $E_P$ (resp. $E_Q$) denotes the special divisor of $\widetilde{X}$ above the point $P=(1,0)$ (resp. $Q=(0,1)$) of $X$. Denote the coefficient of $E_P$ (resp. $E_Q$) in the Weil divisor of $F$ by $\ord_{E_P}(F)$ (resp. $\ord_{E_Q}(F)$). We refer to these as the orders (of vanishing) of $F$ along $E_P$ and $E_Q$. For referencing purposes, we record the following explicit characterization of these orders, which is easily seen through computations in charts:
\begin{lemma}\label{lem: 1st char. of order}
The order of $F$ along $E_P$ is the largest integer $s$ such that $(x-1)^s$ divides $F(x,u(x-1))$ in the polynomial ring $\Q[x,u]$. The order of $F$ along $E_Q$ is the largest integer $s$ such that $(y-1)^s$ divides $F(u(y-1),y)$ in the polynomial ring $\Q[u,y]$.
\end{lemma}
The following definition is inspired by Definition 3.5 of \cite{Dupont}.
\begin{defn}\label{def: integrability}
We say a meromorphic 2-form on $X$ is \emph{integrable} if it is of the form
\[
\omega_{F,t}=  \frac{F dxdy}{g^{t+1}}
\]
for some $F\in\Q[x,y]$ and $t\in\Z_{\geq 0}$ with $\ord_{E_P}(F)\geq t$ and $\ord_{E_Q}(F)\geq t$.
\end{defn}
Note that since the order of $g$ along both $E_P$ and $E_Q$ is 1, the definition above indeed depends on the differential form, and not on the choice of $F$ and $t$. The set of all integrable 2-forms on $X$ is a subspace of the rational vector space of meromorphic 2-forms on $X$, and it contains the space of regular 2-forms on $X$. The terminology \emph{integrable} here refers to integrability over the simplex $\D$, in that the integral of an integrable differential form (in the above sense) converges over $\D$. Indeed, one has the following:
\begin{lemma}\label{lem: integrability condition}
The differential form $\omega_{F,t}$ is integrable if and only if, considered as a meromorphic differential form on $\widetilde{X}$, it is regular on $\widetilde{X}\setminus \widetilde{A}$. (In particular, the integral of an integrable differential form over $\D$ converges.)
\end{lemma}
The proof of this lemma is done by calculations in charts similar to the proof of Lemma \ref{deRham} and is omitted.

We can now formulate a more natural version of Theorem \ref{thm: linear forms}:

\begin{thm}\label{thm: linear forms, ver2}
Let $t$ be an integer $\geq 0$, and $F$ a $\sigma$-invariant element of $\Q[x,y]$ such that $\omega_{F,t}$ is integrable. Then the integral of $\omega_{F,t}$ over $\D$ converges to a $\Q$-linear combination of 1 and $\bG$.    
\end{thm}

This statement is (seemingly, see Remark \ref{rem: ver 1 implies ver 2}) more general than Theorem \ref{thm: linear forms}. Indeed, in view of Lemma \ref{lem: 1st char. of order}, $\omega_{F,t}$ is integrable if $x^ty^t\mid F$, so that condition (ii) of Theorem \ref{thm: linear forms} implies integrability. On the other hand, the condition $x^ty^t\mid F$ is not necessary for integrability of $\omega_{F,t}$. For instance, the differential form 
\[\frac{(1-x-y)(1-x+y)(1+x-y)(1+x+y)}{g^3} \, dx dy\]
is integrable and hence Theorem \ref{thm: linear forms, ver2} asserts that its integral over $\D$ is in $\langle 1,\bG\rangle_\Q$.

\begin{rem}\label{rem: ver 1 implies ver 2}
Despite its more general appearance, Theorem \ref{thm: linear forms, ver2} can be deduced from Theorem \ref{thm: linear forms}. This can be seen easily using the considerations of \S \ref{sec: computation in the general integrable case}.
\end{rem}

\subsection{Symmetric differential forms and the motivic proof of Theorem \ref{thm: linear forms, ver2}}\label{sec: sym diff forms}
Let $F\in \Q[x,y]$ and $t\in\Z$ be such that $\omega_{F,t}$ is integrable. Thanks to its regularity on $\widetilde{X}\setminus\widetilde{A}$ (see Lemma \ref{lem: integrability condition}), $\omega_{F,t}$ defines an element of $M_\textrm{dR}$, which we also denote by $\omega_{F,t}$.
\medskip\par 
Recall from \S \ref{Catalanmotive} that we have an automorphism $\sigma$ of $X$ defined by the formula $(x,y)\mapsto (-y,x)$ at the level of points, which induces the automorphism $\sigma$ of $H^1(A,A\cap B)$ that we used to define the motive $C$ of the Catalan constant. More generally, we will use the notation $\sigma$ to refer to any automorphism of any object (when applicable) that is induced by the automorphism
\[\sigma \colon \Q[x,y]\rightarrow \Q[x,y] \hspace{.5in} (x,y) \mapsto (-y,x).\]
Thus for example, we have automorphisms $\sigma$ of $X$, $A$, $A\cap B$ (hence $X\setminus A$ and the pair $(A,A\cap B)$), $\wX\setminus\wA$, and the motives and the spaces of differential forms on said objects.

\medskip\par 
Let $\Q[x,y]^\sigma$ be the subring of $\Q[x,y]$ consisting of all the $\sigma$-invariant elements:
\[{\Q[x,y]}^\sigma:=\{F\in \Q[x,y]: \sigma F = F\}.\]
Note that for any $t$, the polynomial $F\in \Q[x,y]$ is $\sigma$-invariant if and only if $\omega_{F,t}$ is a $\sigma$-invariant differential form.
\medskip\par 
The first assertion of the following proposition might seem obvious at first glance. However, this is not the case, as the action of $\sigma$ on $H^1(A,A\cap B)$ does not seem to extend to an action on the motive $M$ (as importantly, $\sigma$ does not act on $B$).
\begin{prop}\label{prop: proof of lin form thm, t=0 case}
Let $t\in\Z$. Suppose $F$ is a polynomial in $\Q[x,y]^\sigma$ such that $\omega_{F,t}$ is integrable. Then the element $\omega_{F,t}\in M_\textup{dR}$ belongs to $C_\textup{dR}$. In particular, the statement of Theorem \ref{thm: linear forms, ver2} (and hence Theorem \ref{thm: linear forms}) holds.
\end{prop}

\begin{proof}
First, note that if the element $\omega_{F,t}\in M_\textup{dR}$ belongs to $C_\textup{dR}=\langle \eta,\omega_0\rangle_\Q$, then pairing $\omega_{F,t}$ with the element $\widetilde{\Delta}$ of the Betti homology $M^\vee_{\textup{sing}}$, in view of Remark \ref{rem: regularity and convergence} and the period matrix of $M$ in Proposition \ref{prop: period matrix} we obtain Theorem \ref{thm: linear forms, ver2}. It thus remains to prove that $\omega_{F,t}$ is in $C_\dR$. By definition of $C$, to show this it suffices to show that the element $\res_\wA(\omega_{F,t})$ of $H^1_\dR(A,A\cap B)$ is fixed by $\sigma$.

We will abbreviate $\omega_{F,t}$ by $\omega$ in the rest of the argument. Recall that $\pi: \wX\rightarrow X$ is the blow-up map at $P$ and $Q$ with special divisor $E$. Let $\dbtilde{X}$ be the blow-up of $\wX$ at $\pi^{-1}(-P)$ and $\pi^{-1}(-Q)$, which is the same as the blow-up of $X$ at $\pm P$ and $\pm Q$. Concretely, we take $\dbtilde{X}$ to be the subvariety of $X\times (\mathbb{P}^1)^4$ defined by the appropriate equations (similar to our concrete description of $\wX$, see equation \eqref{p1eq12} and the equations immediately after it). Let $\pi':\dbtilde{X}\rightarrow \wX$ be the blow-up map. Denote the divisor in $\dbtilde{X}$ above $P,Q$ (or $\pi^{-1}(P)$ and $\pi^{-1}(Q)$) by $E_+$ and the divisor in $\dbtilde{X}$ above $-P, -Q$ by $E_-$. The map $\pi'$ maps $E_+$ isomorphically to $E_1:=E$. Set $E_2=E_+\cup E_-$.

Let $B'\subset X$ be the $\sigma$-orbit of $B$ (i.e., the union of the two lines given by $x=0$ and $y=0$, with the four lines given by $1\pm x\pm y =0$). Then $\sigma$ acts on the triple $(X, A, B')$. We denote by $\dbtilde{A}$, $\dbtilde{B}$, and $\dbtilde{B'}$ respectively the strict transforms of $A$, $B$, and $B'$ in $\dbtilde{X}$. The action of $\sigma$ on $(X, A, B')$ extends to an action on $(\dbtilde{X}, \dbtilde{A}, \dbtilde{B'}\cup E_2)$ and hence on the motive $H^2(\dbtilde{X}\setminus \dbtilde{A}, (\dbtilde{B'}\cup E_2)\setminus \dbtilde{A})$. (Note that $\sigma$ acts on $\dbtilde{X}$, but not on $\wX$.)

The differential form $\omega$ is regular on $\wX\setminus \wA$, and hence $\omega$ (or rather, its pullback under $\pi'$, which is given by the same formula as $\omega$ in our model of $\dbtilde{X}$ and is $\sigma$-invariant) is regular on $\dbtilde{X}\setminus (\dbtilde{A}\cup E_-)$. Thanks to its $\sigma$-invariance, $\omega$ is regular on $\dbtilde{X}\setminus \sigma^2(\dbtilde{A}\cup E_-)=\dbtilde{X}\setminus (\dbtilde{A}\cup E_+)$, and hence it is regular on $\dbtilde{X}\setminus \dbtilde{A}$. Just in this argument, let $c_\omega$, $\tilde{c}_\omega$, and $\tilde{c}'_\omega$ be respectively the elements of $M_\dR$, $H^2_\dR(\dbtilde{X}\setminus \dbtilde{A}, (\dbtilde{B}\cup E_2)\setminus \dbtilde{A})$, and $H^2_\dR(\dbtilde{X}\setminus \dbtilde{A}, (\dbtilde{B'}\cup E_2)\setminus \dbtilde{A})$ defined by $\omega$. (One can see similarly to Remark \ref{rem: Xtilde-Atilde is affine} that $\dbtilde{X}\setminus \dbtilde{A}$ is affine. Also, note that $\dbtilde{B'}\cup E_2$ is a simple normal crossing divisor in $\dbtilde{X}$.)

The argument now proceeds in two steps. 

\underline{Step one}: We have a commutative diagram
\[
\begin{tikzcd}[row sep = small]
   H^2_\dR(\dbtilde{X}\setminus \dbtilde{A}, (\dbtilde{B'}\cup E_2)\setminus \dbtilde{A}) \ar[dd] \ar[r, "\res_{\dbtilde{A}}"] & H^1_\dR(\dbtilde{A},\dbtilde{A}\cap (\dbtilde{B'}\cup E_2)) \quad \ar[dd, equal] & \\
   & & H^1_\dR(A,A\cap B), \ar[dl, , equal] \ar[ul, , equal] \\
   H^2_\dR(\dbtilde{X}\setminus \dbtilde{A}, (\dbtilde{B}\cup E_2)\setminus \dbtilde{A})  \ar[r, "\res_{\dbtilde{A}}"] & H^1_\dR(\dbtilde{A},\dbtilde{A}\cap (\dbtilde{B}\cup E_2)) &
\end{tikzcd}
\]
where the vertical map on the left is induced by the identity map on $\dbtilde{X}\setminus \dbtilde{A}$, and the slanted identifications with $H^1(A,A\cap B)$ are induced by $\pi\circ\pi'$ (which gives an isomorphism $\dbtilde{A}\rightarrow A$). There is a natural action of $\sigma$ on every cohomology group in the diagram except $H^2_\dR(\dbtilde{X}\setminus \dbtilde{A}, (\dbtilde{B}\cup E_2)\setminus \dbtilde{A})$, and the identifications in the diagram as well as the top residue map are $\sigma$-equivariant. Since $\omega$ is a $\sigma$-invariant regular form on $\dbtilde{X}\setminus\dbtilde{A}$, the element $\tilde{c}'_\omega$ of $H^2_\dR(\dbtilde{X}\setminus \dbtilde{A}, (\dbtilde{B'}\cup E_2)\setminus \dbtilde{A})$ is $\sigma$-invariant. By the $\sigma$-equivariance of the top residue map, $\res_{\dbtilde{A}}(\tilde{c}'_\omega)$ is $\sigma$-invariant. The vertical map on the left sends $\tilde{c}'_\omega$ to $\tilde{c}_\omega$, so by the commutativity of the diagram we get
\begin{equation}\label{eqp2eq23}
    \res_{\dbtilde{A}}(\tilde{c}_\omega)=\res_{\dbtilde{A}}(\tilde{c}'_\omega) \in H^1_\dR(A,A\cap B)^\sigma.
\end{equation}

\underline{Step two}: We now relate $\res_{\dbtilde{A}}(\tilde{c}_\omega)$ and $\res_{\wA}(c_\omega)$ in $H^1_\dR(A,A\cap B)$. This is done via the following diagram:
\begin{equation}\label{eq motivic proof of main thm}
\begin{tikzcd}
   H^2_\dR(\dbtilde{X}\setminus \dbtilde{A}, (\dbtilde{B}\cup E_2)\setminus \dbtilde{A}) \ar[d, "\iota^\ast" '] \ar[r, "\res_{\dbtilde{A}}"] & H^1_\dR(A,A\cap B) \ar[dd, equal]   \\
  H^2_\dR(\dbtilde{X}\setminus (\dbtilde{A}\cup E_-), (\dbtilde{B}\cup E_2)\setminus (\dbtilde{A}\cup E_-)) &  \\
   M_\dR=H^2_\dR(\wX\setminus \wA, (\wB\cup E)\setminus \wA) \ar[u, equal, "\pi'^\ast"] \ar[r, "\res_{\wA}"] & H^1_\dR(A,A\cap B)
\end{tikzcd}
\end{equation}
The identification marked as $\pi'^\ast$ is induced by the blow-up map $\dbtilde{X}\xrightarrow{\pi '} \wX$, which gives an isomorphism of the relevant two pairs. The map marked as $\iota^\ast$, which is induced by the inclusion map $\dbtilde{X}\setminus (\dbtilde{A}\cup E_-) \xrightarrow{\iota} \dbtilde{X}\setminus \dbtilde{A}$, is also an isomorphism. Indeed, this can be checked in the Betti realization (note that all the maps above and in particular, $\iota^\ast$ are motivic) as a consequence of excision, by a similar argument to the one given in the proof of Lemma \ref{excision lemma}. It follows that $H^2_\dR(\dbtilde{X}\setminus \dbtilde{A}, (\dbtilde{B}\cup E_2)\setminus \dbtilde{A})$ has a basis given by the classes of differential forms $\eta$ and $\omega_j$ ($0\leq j\leq 3$) defined earlier, now considered as differential forms on $\dbtilde{X}$, regular on $\dbtilde{X}\setminus \dbtilde{A}$ and with at most logarithmic singularities along $\dbtilde{A}$. We can now check by a direct computation of the Poincar\'{e} residues (similar to the proof of Lemma \ref{lem: Res of omega_j}) that the images of $\eta$ and $\omega_j$ ($0\leq j\leq 3$) under $\res_{\dbtilde{A}}$ and $\res_{\wA}$ coincide in $H^1_\dR(A,A\cap B)$. It follows that \eqref{eq motivic proof of main thm} is commutative. 

We can now finish the proof. Since $\iota^\ast$ maps $\tilde{c}_\omega$ to $c_\omega$, we have $\res_\wA(c_\omega)=\res_{\dbtilde{A}}(\tilde{c}_\omega)$. Combining with \eqref{eqp2eq23} we get $\sigma$-invariance of $\res_\wA(c_\omega)$.
\end{proof}

\begin{rem}
The case $t=0$ of Proposition \ref{prop: proof of lin form thm, t=0 case} (where the differential form $\omega_{F,t}$ has at most a logarithmic singularity along $\wA$) can be checked by an explicit computation of $\res_\wA(\omega_{F,0})\in H^1(A,A\cap B)$ using the Poincar\'{e} residue of the differential form $\omega_{F,0}$ along $\wA$. In the case $t>0$ however, an explicit computation of $\res_\wA(\omega_{F,t})\in H^1(A,A\cap B)$ is not straightforward, as the Poincar\'{e} residue map for differential forms is only defined for differential forms with at most logarithmic singularity. The proof of Proposition \ref{prop: proof of lin form thm, t=0 case} given above, which avoids residue calculations in the non-logarithmic case, was suggested to us by an anonymous referee. We sincerely thank the referee for this.
\end{rem}

\begin{rem}
One may like to try to deduce the commutativity of \eqref{eq motivic proof of main thm} using naturality of the Gysin long exact sequence, by adding a natural residue map out of
\begin{equation}\label{p2eq24}
H^2_\dR(\dbtilde{X}\setminus (\dbtilde{A}\cup E_-), (\dbtilde{B}\cup E_2)\setminus (\dbtilde{A}\cup E_-))
\end{equation}
to the second row of the diagram. However, one needs to be careful as the target of the natural residue map out of \eqref{p2eq24} is not $H^1_\dR(A,A\cap B)$. Indeed, the triple $(\dbtilde{X}\setminus E_-, (\dbtilde{B}\cup E_2)\setminus E_-, \dbtilde{A}\setminus E_-)$ of the smooth variety $\dbtilde{X}\setminus E_-$, the closed subvariety $(\dbtilde{B}\cup E_2)\setminus E_-$, and the smooth closed subvariety $\dbtilde{A}\setminus E_-$ of pure codimension 1 gives a Gysin long exact sequence with a connecting morphism
\[
\begin{tikzcd}[row sep = small]
H^2_\dR(\dbtilde{X}\setminus (\dbtilde{A}\cup E_-), (\dbtilde{B}\cup E_2)\setminus (\dbtilde{A}\cup E_-))  \ar[r, "\res_{\dbtilde{A}\setminus E_-}"] & H^1_\dR(\dbtilde{A}\setminus E_-, (\dbtilde{A}\cap E_2)\setminus E_-) \ar[d,equal] \\
& H^1_\dR(A\setminus\{-P,-Q\}, \{P,Q\})
\end{tikzcd}
\]
in the de Rham realization (note that $\dbtilde{A}$ does not meet $\dbtilde{B}$). Notice that the codomain is different from the codomains of $\res_{\dbtilde{A}}$ and $\res_\wA$ in \eqref{eq motivic proof of main thm}. One could try to proceed by adding this map to the diagram \eqref{eq motivic proof of main thm}. We also have an obvious map from $H^1_\dR(A,A\cap B)$ to $H^1_\dR(A\setminus\{-P,-Q\}, \{P,Q\})$. We obtain a new diagram
\[
\begin{tikzcd}
   H^2_\dR(\dbtilde{X}\setminus \dbtilde{A}, (\dbtilde{B}\cup E_2)\setminus \dbtilde{A}) \ar[d, "\iota^\ast" '] \ar[r, "\res_{\dbtilde{A}}"] & H^1_\dR(A,A\cap B) \ar[d]   \\
  H^2_\dR(\dbtilde{X}\setminus (\dbtilde{A}\cup E_-), (\dbtilde{B}\cup E_2)\setminus (\dbtilde{A}\cup E_-)) \ar[r, "\res_{\dbtilde{A}\setminus E_-}"] &  H^1_\dR(A\setminus\{-P,-Q\}, \{P,Q\}), \\
   M_\dR=H^2_\dR(\wX\setminus \wA, (\wB\cup E)\setminus \wA) \ar[u, equal, "\pi'^\ast"] \ar[r, "\res_{\wA}"] & H^1_\dR(A,A\cap B) \ar[u]
\end{tikzcd}
\]
which commutes for functoriality reasons. However, this does not imply commutativity of the original diagram as the map $H^1_\dR(A,A\cap B)\rightarrow H^1_\dR(A\setminus\{-P,-Q\}, \{P,Q\})$ has a nontrivial kernel. (In fact, it is not difficult to compute this kernel: it is generated by $\ow_1-\ow_2$ and $4\ow_3-(\ow_1+\ow_2)$; notice that the kernel is not contained in the $\sigma$-invariant subspace.)
\end{rem}

\subsection{Further comments on symmetric differential forms}\label{sec: sym diff forms further comments}
In this subsection we make some preliminary observations that will prepare us for \S \ref{sec: lin forms II} and \S \ref{sec: lin forms III}.

\medskip\par 
For every integers $m,n$, consider the following elements of $\Q(x,y)$:
\begin{align}\label{eq: def of phi, psi}
\begin{split}
\phi_{m,n} & := x^my^n+x^ny^m \\
\psi_{m,n} & := x^my^n-x^ny^m. 
\end{split}
\end{align}
The proof of the following lemma is clear and is omitted.

\begin{lemma}\label{lem: basis for sigma-invariants} \ 
\begin{itemize}
\item[$(a)$] Let $F=\sum\limits_{m,n}c_{m,n}x^my^n\in \Q[x,y]$. Then $\sigma F = F$ if and only if for every nonnegative integers $m,n$, we have $c_{n,m}=(-1)^mc_{m,n}$. In particular, if $F$ is $\sigma$-invariant and $c_{m,n}\neq 0$, then $m\equiv n \pmod 2$.
\item[$(b)$] The collection of polynomials 
    \[
    \{\phi_{m,n}: m\geq n\geq 0;\, 2\mid m;\, 2\mid n\} \ \cup \ \{\psi_{m,n}: m>n>0;\, 2\nmid mn\}
    \]
forms a basis of ${\Q[x,y]}^\sigma$.
\end{itemize}
\end{lemma}

Recall from Remark \ref{rem: automorphism tau} that the involution $\tau$ of $X$ defined by $(x,y)\mapsto (y,x)$ extends to an involution on $\wX$, which then restricts to an automorphism of the pair $(\wX\setminus \wA, (\wB\cup E)\setminus \wA)$ and hence induces an automorphism of the motive $M$. With abuse of notation, we denote all of these, as well as the corresponding maps at the level of coordinate rings and differential forms on $X$ and $\wX$, by $\tau$. In particular, the map
\[
\tau \colon \Q[x,y]\rightarrow \Q[x,y]
\]
is given by
\[
(x,y)\mapsto (y,x).
\]
We have a decomposition into eigenspaces
\[
M= M^\tau\oplus M^{\tau=-1},
\]
and $C\subset M^{\tau=-1}$ (see Remark \ref{rem: automorphism tau}). There is also a decomposition
\begin{equation}\label{p2eq22}
\Q[x,y] = {\Q[x,y]}^\tau \oplus {\Q[x,y]}^{\tau=-1}
\end{equation}
into eigenspaces. For every $F\in \Q[x,y]$ and every integer $t\geq 0$, we have
\begin{equation}\label{p2eq1}
\tau \omega_{F,t} = (\tau F) \, \frac{dydx}{g^{t+1}} = \omega_{-\tau F, t}
\end{equation}
(see \eqref{p2eq11} to recall the definition of $\omega_{F,t}$). These equalities are equalities of differential forms on $\wX$ (and on $X$). Since $\tau$ acts on the pair $(\wX\setminus \wA, (\wB\cup E)\setminus \wA)$, assuming $\omega_{F,t}$ is integrable (or equivalently, regular on $\wX\setminus\wA$, see Lemma \ref{lem: integrability condition}; note that this condition holds automatically when $t=0$), the equalities descend to equalities in $M_\textrm{dR}$. From this we immediately have:
\begin{lemma}\label{lem: tau-anti inv F}
Suppose $F\in {\Q[x,y]}^{\tau=-1}$ and $\omega_{F,t}$ is integrable. Then $\omega_{F,t}\in M_\textup{dR}^\tau$.
\end{lemma}

\begin{rem}\label{rem: interpretation of tau-anti invariance lem}
Lemma \ref{lem: tau-anti inv F} is a geometric upgrading of the fact, observable through a trivial elementary computation, that if $F\in {\Q[x,y]}^{\tau=-1}$, then the integral of $\omega_{F,t}$ over $\D$ vanishes provided that it convergences. Indeed, recall from Remark \ref{rem: automorphism tau} that $M_\textrm{dR}^\tau$ is generated by $\omega_1-\omega_2$. From the period matrix of $M$ given in Proposition \ref{prop: period matrix}, the integral of $\omega_1-\omega_2$ over $\D$ is zero.
\end{rem}

Note that
\begin{align*}
&{\Q[x,y]}^\tau = \langle \phi_{m,n}: m\geq n\geq 0\rangle_\Q\\
&{\Q[x,y]}^{\tau=-1} = \langle \psi_{m,n}: m> n\geq 0\rangle_\Q
\end{align*}
(where $\phi_{m,n}$ and $\psi_{m,n}$ are as in \eqref{eq: def of phi, psi}). Combining with Lemma \ref{lem: basis for sigma-invariants}, we obtain a decomposition
\begin{equation}\label{p2eq2}
{\Q[x,y]}^\sigma = {\Q[x,y]}^{\sigma,\tau} \oplus {\Q[x,y]}^{\substack{\sigma=1 \\ \tau=-1}}
\end{equation}
with
\begin{equation}\label{p2eq21}
{\Q[x,y]}^{\sigma,\tau}:=\{F\in \Q[x,y]: \sigma F=F=\tau F\} = \langle \phi_{m,n}: m\geq n\geq 0; 2\mid m; 2\mid n\rangle_\Q
\end{equation}
and
\[
{\Q[x,y]}^{\substack{\sigma=1 \\ \tau=-1}}:=\{F\in \Q[x,y]: \sigma F=F=-\tau F\} = \langle \psi_{m,n}: m> n> 0; 2\nmid mn \rangle_\Q.
\]
The following result will be useful in \S \ref{sec: lin forms II} and \S \ref{sec: lin forms III}.
\begin{lemma}\label{lem: sigma=1, tau=-1 F}
Let $F\in {\Q[x,y]}^{\substack{\sigma=1 \\ \tau=-1}}$. Suppose $\omega_{F,t}$ is integrable. Then $\omega_{F,t}$ vanishes in $M_\textup{dR}$.
\end{lemma}
\begin{proof}
By Proposition \ref{prop: proof of lin form thm, t=0 case}, $\omega_{F,t}$ is in $C_\textrm{dR}$. On the other hand, by Lemma \ref{lem: tau-anti inv F}, $\omega_{F,t}$ is in $M_\textrm{dR}^\tau$. The intersection of $C_\textrm{dR}$ and $M_\textrm{dR}^\tau$ is zero (see Remark \ref{rem: automorphism tau}).
\end{proof}

We finish this discussion with one last preparation for the explicit computations of the next sections. As before, let $K=\Q(i)$. The actions of $\sigma$ and $\tau$ on $\Q[x,y]$ extend $K$-linearly to automorphisms of $K[x,y]$.

We will identify
\begin{equation}\label{p2eq3}
K[x,y] \cong K[z,\tw] \ \ \ \ \ \ \text{via} \ \ \ \ \ \ z= x+iy \ \ \ \text{and} \ \ \tw=x-iy, \ \ \text{fixing $K$}.
\end{equation}
This identification induces the earlier identification of $A_K=\operatorname{Spec}K[x,y]/(g)$ with $\mathbb{G}_{m, K}= \operatorname{Spec}K[z,\tw]/(1-z\tw)$. The actions of $\sigma$ and $\tau$ on $x,y$ translate under \eqref{p2eq3} to actions on $z,\tw$ given by
\begin{align*}
\sigma \colon \ \ & (z,\tw)\mapsto (iz, -i\tw)\\
\tau \colon \ \ & (z,\tw)\mapsto (i\tw, -iz).
\end{align*}
Moreover, the action of complex conjugation on $K[x,y]$ translates to an action on $K[z,\tw]$ which conjugates the coefficients and interchanges $z$ and $\tw$. From these we obtain the following characterization of $\{\sigma,\tau\}$-invariance in terms of the variables $z$ and $\tw$.
\begin{lemma}\label{lem: invariance in terms of z,w}
Let $F= \sum\limits_{k,l} \lambda_{k,l} z^k\tw^l\in K[z,\tw]$. Then we have:
\begin{enumerate}[label=$(\alph*)$] 
\item $F$ is fixed by $\sigma$ \ \ $\Longleftrightarrow$ \ \ If $\lambda_{k,l}\neq 0$, then $k\equiv l \pmod 4$.
\item $F$ is fixed by $\tau$ \ \ $\Longleftrightarrow$ \ \ $\lambda_{l,k}=i^{k-l}\lambda_{k,l}$ for all $k,l$.
\item $\tau F=-F$ \ \ $\Longleftrightarrow$ \ \ $-\lambda_{l,k}=i^{k-l}\lambda_{k,l}$ for all $k,l$.
\item $F\in \Q[x,y]$ \ \ $\Longleftrightarrow$ \ \ $\lambda_{l,k}=\overline{\lambda_{k,l}}$ for all $k,l$.
\end{enumerate}
\end{lemma}

\section{Linear forms in 1 and $\bG$ -- II: Calculation of the coefficients in the simple pole case}\label{sec: lin forms II}

The objective of the rest of the paper is to make Theorem \ref{thm: linear forms, ver2} constructive; that is, to explicitly calculate the coefficients of 1 and $\bG$ in the linear forms of Theorem \ref{thm: linear forms, ver2}. In this section, we will deal with the case $t=0$, which will turn out to be the core case. The case of the theorem for higher order poles will be handled in \S \ref{sec: lin forms III}.

Given any $F\in\Q[x,y]$, for simplicity, we shall write $\omega_F$ for $\omega_{F,0}=Fdxdy/g$. Recall from Proposition \ref{prop: proof of lin form thm, t=0 case} that if $F\in\Q[x,y]$ is $\sigma$-invariant, then $\omega_F$ represents an element of $C_\textrm{dR}$ (as the integrability condition holds for every $F$ when $t=0$). In this section, we shall explicitly calculate the rational numbers $a_F, b_F$ such that
\[
\omega_F = a_F \, \eta + b_F \, \omega_0
\]
in $C_\textrm{dR}$. In particular, this will make the case $t=0$ of Theorem \ref{thm: linear forms, ver2} constructive.

We will deal with the computation of $b_F$ and $a_F$ separately. The computation of $b_F$ is easier and will be handled first.

\subsection{Calculation of the coefficient of $\omega_0$}
This amounts to a Poincar\'{e} residue calculation. We will do this for any polynomial $F\in\Q[x,y]$, with no $\sigma$-invariance assumption. Throughout, we fix an arbitrary $F\in \Q[x,y]$. Then $\omega_F$ defines an element of $M_\textrm{dR}$ and hence it can uniquely be written as
\[
\omega_F \ = \ a_F\, \eta + \, \text{(a $\Q$-linear combination of $\omega_1,\omega_2,\omega_3$)} \, + \, b_F \, \omega_0
\]
for some $a_F$ and $b_F$ in $\Q$.

There is a commutative diagram
\[
\begin{tikzcd}
    M_\textrm{dR} \ar[r, "\res_\wA"] \ar[d] & H^1_\textrm{dR}(A,A\cap B) \ar[d] \\
    H^2_\textrm{dR}(\wX\setminus \wA) \ar[r, "\res_\wA"] &  H^1_\textrm{dR}(A),
\end{tikzcd}
\]
where the vertical maps are the natural maps for relative cohomology, and we have used the identifications \eqref{isom1} and $H^1_{\mathrm{dR}}(\wA)\cong H^1_{\mathrm{dR}}(A)$ given by the blow-up. The top residue map sends $\eta$ to zero and $\omega_j$ to $\overline{\omega}_j$ for $0\leq j\leq 3$ (see Lemma \ref{lem: Res of omega_j}). The vertical map on the right sends $\overline{\omega}_j$ to zero for $1\leq j\leq 3$ and sends $\overline{\omega}_0=dx/(2y)$ to its namesake in $H^1_\textrm{dR}(\wA)$ (see Remark \ref{rem: wt filt of H^1_dR(A, A cap B)}). Passing $\omega_F$ along the left vertical map to consider it as an element of $H^2_\textrm{dR}(\wX\setminus\wA)$, we thus have
\begin{equation}\label{p2eq4}
\res_\wA(\omega_F) = b_F \overline{\omega}_0 \ \ \ \text{in} \ \ H^1_\textrm{dR}(A) = \langle \overline{\omega}_0\rangle_\Q. 
\end{equation}
To find $b_F$, we compute $\res_\wA(\omega_F)$ in $H^1_\textrm{dR}(A)$. As a differential form on $\wX$ with at most a logarithmic singularity along $\wA$,
\[
\omega_F = \frac{Fdxdy}{g} = \frac{dg}{g}\wedge \frac{Fdx}{2y},
\]
so that
\[
\res_\wA(\omega_F) = \frac{Fdx}{2y}
\]
as 1-forms on $A$ ($\cong \wA$). To compute its cohomology class in $H^1_\textrm{dR}(A)$, we first extend the base field to $K=\Q(i)$. Then
\[A_K \cong \mathbb{G}_{m,K} = \operatorname{Spec}K[z,\tw]/(1-z\tw),\] 
the isomorphism given by \eqref{p2eq3}. The differential form $Fdx/(2y)$ on $A$ becomes
\[
F(\frac{z+\tw}{2}, \frac{z-\tw}{2i})\bigg\vert_{\tw=z^{-1}} \cdot \frac{idz}{2z}
\]
on $\mathbb{G}_{m, K}$. It follows that in $H^1_\textrm{dR}(A_K)$ (and hence in $H^1_\textrm{dR}(A)$),
\[
\res_\wA(\omega_F) \ = \ \begin{pmatrix} \text{  the constant term of} \\
\  \displaystyle{F(\frac{z+\tw}{2}, \frac{z-\tw}{2i})\bigg\vert_{\tw=z^{-1}}} \
\end{pmatrix} \cdot \, \frac{idz}{2z} \ = \begin{pmatrix} \text{ the constant term of} \\
\  \displaystyle{F(\frac{z+\tw}{2}, \frac{z-\tw}{2i})\bigg\vert_{\tw=z^{-1}}} \
\end{pmatrix} \cdot \, \overline{\omega}_0.
\]
Combining with \eqref{p2eq4}, we obtain the following:
\begin{prop}\label{prop: computation of b_F}
Let $F\in \Q[x,y]$. Writing
\[
\omega_F \ = \ a_F\, \eta + \, (\text{a $\Q$-linear combination of $\omega_1,\omega_2,\omega_3$}) \, + \, b_F \, \omega_0
\]
in $M_\textrm{dR}$, we have
\[
b_F \ = \ \text{the constant term of} \ F(\frac{z+\tw}{2}, \frac{z-\tw}{2i})\bigg\vert_{\tw=z^{-1}} \, .
\]
\end{prop}

\subsection{Calculation of the coefficient of $\eta$}\label{sec: coefficient of eta, log case}
We now assume $F\in\Q[x,y]$ is $\sigma$-invariant, so that we have $\omega_F = a_F\eta+b_F\omega_0$ in $C_\textrm{dR}$ for some unique $a_F, b_F\in\Q$. Our goal is to compute $a_F$. 

Since we have already computed $b_F$ in the previous subsection, we can try to compute $a_F$ by pairing $\omega_F\in C_\textrm{dR}$ with the element $\widetilde{\Delta}$ of $C^\vee_{\textup{sing}}$, which amounts to integration on $\Delta$. We obtain
\begin{equation}\label{p2eq5}
\int_\Delta (\omega_F-b_F \, \omega_0) = a_F\int_\Delta \eta  \ = a_F.
\end{equation}
However, the task is far from over since we need to evaluate the integral on the left explicitly. This calculation is done in this subsection.

Thanks to the decomposition \eqref{p2eq2} and Lemma \ref{lem: sigma=1, tau=-1 F}, we may, and indeed shall, assume that
\[
F\in \Q[x,y]^{\sigma, \tau}.
\]
Start by writing
\begin{equation}\label{p2eq9}
F(\frac{z+\tw}{2},\frac{z-\tw}{2i}) \ = \sum\limits_{k,l} \, \lambda_{k,l} \,z^k\tw^l \, \in \, K[z,\tw] 
\end{equation}
(as before, $K=\Q(i)$). By Lemma \ref{lem: invariance in terms of z,w}, we have
\[
\begin{cases} \ \text{If $\lambda_{k,l}\neq 0$, then $k\equiv l \pmod 4$, and} \\
\ \lambda_{k,l}=\lambda_{l,k}\in \Q \ \ \text{for all $k$ and $l$}.  
\end{cases}
\]
By Proposition \ref{prop: computation of b_F},
\[
b_F = \sum\limits_k \, \lambda_{k,k}.
\]
Putting these together, we have
\[
F-b_F \ = \sum\limits_{\substack{k>l\geq 0 \\ k\equiv l \pmod4}} \lambda_{k,l}\, (z^k\tw^l+z^l\tw^k)  \ + \ \sum\limits_{k>0} \lambda_{k,k} \, (z^k\tw^k-1) \ \in \ \Q[z,\tw].
\]
Multiplying with $dxdy/g$ and integrating over $\Delta$, in view of \eqref{p2eq5} we get
\begin{equation}\label{p2eq6}
a_F \ = \sum\limits_{\substack{k>l\geq 0 \\ k\equiv l \pmod4}} \lambda_{k,l}\, \int_\Delta \frac{z^k\tw^l+z^l\tw^k}{1-z\tw}dxdy \ + \ \sum\limits_{k>0} \lambda_{k,k} \int_\Delta \frac{z^k\tw^k-1}{1-z\tw}dxdy.
\end{equation}

Before proceeding any further, we note that the integrals in \eqref{p2eq6} all converge to rational numbers. This is obvious for the integrals in the second sum, but the following argument applies to all the integrals: Every 2-form in these integrals is regular on $\wX\setminus \wA$ and gives an element of $M_\textrm{dR}$. The polynomials $z^k\tw^l+z^l\tw^k$ for $k,l$ as in the first sum and $z^k\tw^k-1$ are in $\Q[x,y]^{\sigma,\tau}$ (see Lemma \ref{lem: invariance in terms of z,w}). In particular, by Proposition \ref{prop: proof of lin form thm, t=0 case}, the elements $(z^k\tw^l+z^l\tw^k)dxdy/g$ for $k,l$ as above and $(z^k\tw^k-1)dxdy/g$ of $M_\mathrm{dR}$ live in $C_\textrm{dR}$. Moreover, it follows from Proposition \ref{prop: computation of b_F} that in fact, all these elements belong to $\langle \eta\rangle_\Q\subset C_\textrm{dR}$.

We now proceed to compute the integrals on the right-hand side of \eqref{p2eq6}. Once this is done, $a_F$ can be computed explicitly from $F$.

\subsubsection{Computation of the integrals of the first type}\label{subsubsec: integrals of 1st type}
Fix integers $k,l$ with $k>l$ and $k\equiv l \pmod 4$. Our goal is to compute
\[
\int_\Delta \frac{z^k\tw^l+z^l\tw^k}{1-z\tw}dxdy.
\]
A direct computation shows
\[
\frac{z^k\tw^l+z^l\tw^k}{1-z\tw}dzd\tw = \frac{1}{k-l} \, d\bigm(\frac{z^k\tw^l-z^l\tw^k}{1-zw}\cdot(zdw+wdz) \bigm).
\]
Thus for any 2-chain $\Delta'\subset \Delta$ which avoids $\{P,Q\}$, we can write
\begin{align}
\int_{\Delta'} \frac{z^k\tw^l+z^l\tw^k}{1-z\tw}dxdy & \ \ = \ \frac{i}{2}\int_{\Delta'}\frac{z^k\tw^l+z^l\tw^k}{1-z\tw}dzd\tw \notag\\
& \stackrel{\text{Stokes}}{=} \frac{i}{2(k-l)} \int_{\partial \Delta'} \frac{z^k\tw^l-z^l\tw^k}{1-zw}\cdot(zdw+wdz) \notag\\
& \ \ = \frac{i}{k-l} \int_{\partial \Delta'}\frac{(x+iy)^k(x-iy)^l-(x+iy)^l(x-iy)^k}{1-x^2-y^2}\cdot (xdx+ydy) \,.\label{p2eq13}
\end{align}
The differential forms on the two sides are both regular on $\wX\setminus\wA$; we already know this for the 2-form on the left, and it is easily checked in charts (e.g., using the charts of the proof of Lemma \ref{deRham}) for the 1-form on the right, using the fact that $k$ and $l$ have the same parity. Letting $\Delta'\rightarrow \Delta$, we get
\begin{equation}\label{p2eq18}
\int_{\Delta} \frac{z^k\tw^l+z^l\tw^k}{1-z\tw}dxdy = \frac{i}{k-l} \int_{\partial \wD}\underbrace{\bigm((x+iy)^k(x-iy)^l-(x+iy)^l(x-iy)^k\bigm)}_{f} \, \frac{xdx+ydy}{1-x^2-y^2}.
\end{equation}
(Explanation for the integral on the right: use functoriality to write the integral on the right hand side of \eqref{p2eq13} as an integral on $\pi^{-1}(\partial \D')$ (where $\pi$ = the blow-up map), which then tends to the integral on $\partial \wD$ as $\Delta'\rightarrow \Delta$.) 
\begin{figure}
    \centering
    \includegraphics[width=0.5\linewidth]{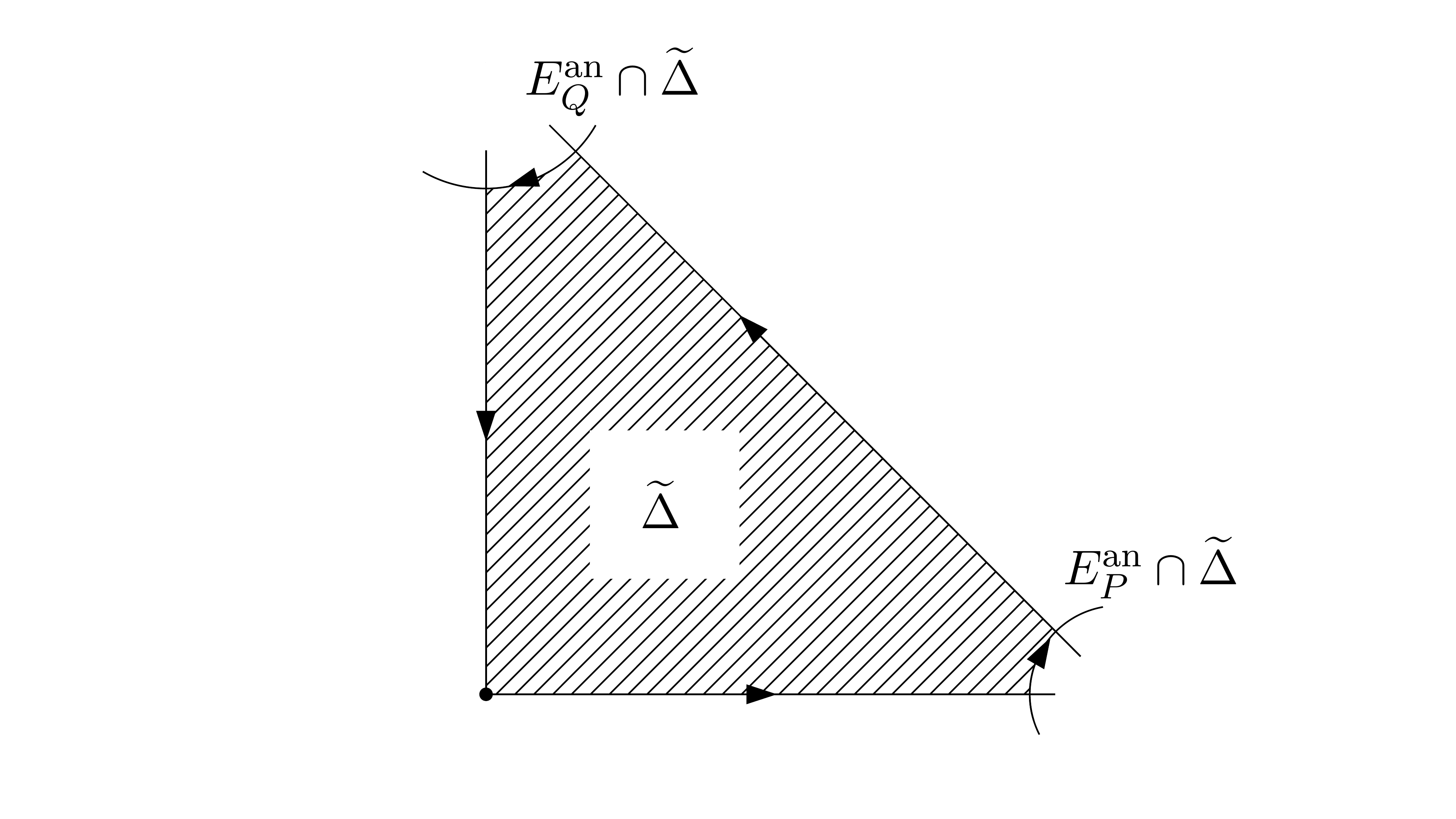}
    \vspace{-.3in}
    \caption{The boundary of $\wD$}
    \label{fig: partial wD}
\end{figure}
The integral on the right hand side of \eqref{p2eq18} now breaks into the sum of the integrals of $f\cdot(xdx+ydy)/g$ (where $f$ is the polynomial marked in the previous equation) over $E_P^\an\cap \wD$, $E_Q^\an\cap \wD$ (where $E_P$ and $E_Q$ are the special divisors, as before), and the rest of $\partial\wD$, all oriented appropriately (see Figure \ref{fig: partial wD}). By functoriality, the latter integral is the integral over $\partial \D$. On the other hand, a computation in charts shows that $f\cdot (xdx+ydy)/g$ vanishes on $E_P\setminus \wA$ and $E_Q\setminus \wA$. In conclusion, we get 
\[
\int_{\Delta} \frac{z^k\tw^l+z^l\tw^k}{1-z\tw}dxdy = \frac{i}{k-l} \int_{\partial \D}  f\cdot\, \frac{xdx+ydy}{1-x^2-y^2},
\]
where the integral on the right converges absolutely on each line segment in $\partial \D$. Now note that the differential form on the right vanishes on the vertical and horizontal line segments in $\partial \Delta$. Parametrizing the straight line segment from $P$ to $Q$ by $(1-y,y)$ as $y$ goes from 0 to 1, we have $g=2y(1-y)$ on this line segment, so that the last equation reads
\[
\int_{\Delta} \frac{z^k\tw^l+z^l\tw^k}{1-z\tw}dxdy = \frac{i}{k-l} \int_0^1 f(1-y,y) \, \frac{(2y-1)dy}{2y(1-y)}.
\]
The polynomial $f(x,y)$ has coefficients in $2i\,\Z$. We observe that $f(1-y,y)$ vanishes at $y=0$ and $y=1$ (as $k$ and $l$ have the same parity). It follows that $i\cdot f(1-y,y)/(2y(1-y))$ is a polynomial in $y$ of degree $k+l-2$ and with coefficients in $\Z$. This completes our computation of the integral on the left. It also gives us an upper bound for its denominator. For referencing purposes, we summarize our conclusions:

\begin{lemma}\label{lem: 1st type integrals, log case}
Let $k,l$ be integers with $k>l\geq 0$ and $k\equiv l \pmod 4$. Then 
\begin{equation}\label{p2eq7}
\int_{\Delta} \frac{z^k\tw^l+z^l\tw^k}{1-z\tw}dxdy \ = \ \frac{1}{k-l} \int_0^1 i\cdot f(1-y,y) \, \frac{(2y-1)dy}{2y(1-y)},
\end{equation}
where 
\[
f(x,y) = (x+iy)^k(x-iy)^l-(x+iy)^l(x-iy)^k.
\]
The integrand function $i\cdot f(1-y,y)\cdot (2y-1)/(2y(1-y))$ is a polynomial of degree $k+l-1$ with coefficients in $\Z$. In particular, the integral \eqref{p2eq7} is a rational number whose lowest denominator divides
\[
(k-l) \, L_{k+l},
\]
where (here and in what follows) for any positive integer $m$, the notation $L_m$ means the least common multiple of $1,\ldots,m$.
\end{lemma}

\begin{rem}\label{rem: integrals of 1st type only need same parity k,l}
As the reader may have noticed, the argument in this subsection only needs $k,l$ to have the same parity. However, note that the case $k-l\equiv 2 \pmod 4$ of the integrals of the first type is trivial, as $z^k\tw^l+z^l\tw^k$ is then in $\Q[x,y]^{\tau=-1}$ and hence $\frac{z^k\tw^l+z^l\tw^k}{1-z\tw}dxdy$ is in $M_\textrm{dR}^\tau$, so its integral over $\Delta$ is zero (see Lemma \ref{lem: tau-anti inv F} and Remark \ref{rem: interpretation of tau-anti invariance lem}).
\end{rem}

\subsubsection{Computation of the integrals of the second type}
We now discuss the second type of integral appearing in \eqref{p2eq6}, i.e., an integral of the form
\begin{equation}\label{p2eq8}
\int_\Delta \frac{z^k\tw^k-1}{1-z\tw} \, dxdy
\end{equation}
where $k$ is a positive integer. The integrand function is now a polynomial so the situation is easy. By trivial algebraic manipulations and Stokes' theorem, the integral above is equal to
\begin{align*}
\frac{i}{2} \int_{\partial\Delta} \left(\sum_{j=0}^{k-1} \frac{\tw^{j+1}z^j}{j+1} \right)dz \ = \ \frac{i}{2} \int_{\partial\Delta} \left(\sum_{j=0}^{k-1} \frac{(x-iy)^{j+1}(x+iy)^j}{j+1} \right)(dx+idy).
\end{align*}
From this we see easily that the {\it sum} of the integrals along the vertical and the horizontal line segments in $\partial \Delta$ is zero. Denoting the line segment from $P$ to $Q$ by $[P,Q]$, the last integral reads
\[
\frac{i}{2} \int_{[P,Q]} \left(\sum_{j=0}^{k-1} \frac{(x^2+y^2)^j}{j+1} \right)(x-iy)(dx+idy) \ = \ \frac{1}{2} \int_{[P,Q]} \left(\sum_{j=0}^{k-1} \frac{(x^2+y^2)^j}{j+1} \right)(ydx-xdy),
\]
since the value of the integral \eqref{p2eq8} is real. Parametrizing $[P,Q]$ by $(1-y,y)$ as before, we obtain the following:
\begin{lemma}\label{lem: 2nd type integrals, log case}
Let $k$ be a positive integer. Then
\[
\int_\Delta \frac{z^k\tw^k-1}{1-z\tw} \, dxdy \ = \ \frac{-1}{2} \int_0^1 \ \sum_{j=0}^{k-1} \frac{1}{j+1} \, (2y^2-2y+1)^j \, dy. 
\]
In particular, this integral is a rational number whose lowest denominator divides
\[
2 \, L_k \, L_{2k-1}.
\]
\end{lemma}

\subsection{Summary: Explicit version of Theorem \ref{thm: linear forms} in the simple pole case}
For future referencing and to help the reader navigate through the arguments, we combine the results of the paper up to this point to record a more explicit version of Theorem \ref{thm: linear forms} in the case $t=0$. 

\begin{thm}\label{thm: lin forms thm, t=0, explicit version}
Let $F\in \Z[x,y]$ be $\sigma$-invariant and of (total) degree $N$. Then
\[
\int_\Delta \frac{F \, dxdy}{1-x^2-y^2} \ = \ a_F+b_F\,\bG,
\]
where $a_F$ and $b_F$ are rational numbers such that $2^N \, b_F$ and $2^{N+2}\,L_N \,L_{N/2} \, a_F$ are integers. (Recall that $L_m$ refers to the least common multiple of $1,\ldots, m$. Note that $N$ is even by $\sigma$-invariance.)
\end{thm}

\begin{proof}
By Proposition \ref{prop: proof of lin form thm, t=0 case}, the integral converges to a $\Q$-linear combination of 1 and $\bG$. To compute $a_F$ and $b_F$, by Lemma \ref{lem: sigma=1, tau=-1 F} we may replace $F$ by the $\{\sigma,\tau\}$-invariant component of it, i.e., the part consisting of the monomials $x^my^n$ where $m$ and $n$ are even (see \eqref{p2eq2} and \eqref{p2eq21}, and \eqref{eq: def of phi, psi}). The coefficients $a_F$ and $b_F$ can now be computed explicitly by Proposition \ref{prop: computation of b_F}, Lemma \ref{lem: 1st type integrals, log case} and Lemma \ref{lem: 2nd type integrals, log case} (see also equations \eqref{p2eq9} and \eqref{p2eq6}). The assertion that $2^N\,b_F$ is an integer follows immediately from Proposition \ref{prop: computation of b_F}, in view of the fact that the polynomial
\[F(\frac{z+\tw}{2}, \frac{z-\tw}{2i}) \]
in $z,\tw$ has coefficients in $2^{-N}\,\Z$. The assertion that $2^{N+2}\,L_N \,L_{N/2} \, a_F$ is an integer follows from the results about the lowest denominators of the integrals in Lemmas \ref{lem: 1st type integrals, log case} and \ref{lem: 2nd type integrals, log case}. Indeed, $2L_N L_{N/2}$ is a common denominator for the integrals in Lemma \ref{lem: 2nd type integrals, log case} for the range of $k$ that may appear in the expansion \eqref{p2eq6}. On the other hand, $2^2L_N L_{\lfloor N/4\rfloor}$ is a common denominator for the integrals in Lemma \ref{lem: 1st type integrals, log case} for the values of $k,l$ that may appear. Note that we used the fact that $4\mid k-l$ if the integral of Lemma \ref{lem: 1st type integrals, log case} appears in the expansion \eqref{p2eq6}.

\end{proof}

\medskip\par 

\begin{example}
Let $F=x^2y^2$. Then $F$ is $\sigma$-invariant, so that $\omega_F=a_F\eta+b_F\omega_0$ in $C_\textrm{dR}$. We have
\[
F(\frac{z+\tw}{2}, \frac{z-\tw}{2i}) =\frac{-1}{16}(z^4-2z^2\tw^2+\tw^4).
\]
Thus $b_F=1/8$ by Proposition \ref{prop: computation of b_F}. Note that $F$ is also fixed by $\tau$. Calculating $a_F$ as described above, we get $a_F=-5/48$.
\end{example}

\begin{rem}
We have proved a stronger statement than the one stated in Theorem \ref{thm: lin forms thm, t=0, explicit version}, in that our proof gives the formula $\omega_F = a_F\, \eta +b_F\, \omega_0$ at the level of $C_\textrm{dR}$.
\end{rem}

\begin{rem}\label{rem: chebychev}
The explicit bounds on the denominators of $a_F$ and $b_F$ in Theorem \ref{thm: lin forms thm, t=0, explicit version} would allow us to establish the irrationality of $\bG$ if there exist polynomials with integer coefficients and sufficiently small sup norms $\|~\|_S$ on the square $S=\cup_{i=0}^3\sigma^i\Delta$ (with vertices $\pm P, \pm Q$). More precisely, suppose there exists a sequence of polynomials $(G_j)_{j\in \Z_{>0}}$ in $\Z[x,y]$ with $\deg(G_j)\rightarrow \infty$ as $j\rightarrow \infty$ such that 
\[
(2(e+\epsilon)^{3/2})^{\deg(G_j)} \ \|G_j\|_S \ \rightarrow 0
\]
as $j\rightarrow \infty$ for some fixed $\epsilon>0$. Then the polynomials $F_j=\sum_{i=0}^3 \sigma^i(G_j^2)\in \Z[x,y]$ are non-negative on $\Delta$, $\sigma$-invariant, and satisfy
\begin{equation}\label{p2eq10}
(2(e+\epsilon)^{3/2})^{\deg(F_j)} \ \|F_j\|_S \ \rightarrow 0
\end{equation}
as $j\rightarrow \infty$. Thanks to Theorem \ref{thm: lin forms thm, t=0, explicit version} we will then have
\[
0 \, < \ \int_\Delta \frac{F_j\,dxdy}{1-x^2-y^2} \ = \ a_j+b_j\,\bG,
\]
where $a_j, b_j$ are rational numbers with a common denominator at most
\[
2^{\deg(F_j)+2}L_{\deg(F_j)}L_{\deg(F_j)/2} \ \ll \ (2(e+\epsilon)^{3/2})^{\deg(F_j)} 
\]
by the Prime Number Theorem. Together with \eqref{p2eq10} this would imply irrationality of $\bG$.
\medskip\par 
The question of whether or not a sequence of polynomials $G_j$ as above exists amounts to whether the \emph{bivariate integer Chebyshev constant of $S$} in the sense of \cite[Definition 1.1]{BP} is smaller than $1/(2e^{3/2})$. The bivariate integer Chebyshev constant of $S$ does not seem to be known, but it appears highly unlikely that this constant is as small as needed for the above argument.
\end{rem}

\section{Linear forms in 1 and $\bG$ -- III: Calculation of the coefficients for non-simple poles}\label{sec: lin forms III}
\subsection{Overview}
The goal of this section is to make Theorem \ref{thm: linear forms, ver2} constructive in the case $t>0$, that is, when our differential form
\[
\omega_{F,t}=  \frac{F dxdy}{g^{t+1}}
\]
may have a non-simple pole along $A$. The motivation for this is the hope that it may be easier to find better-behaved linear forms by allowing the order of the pole to grow.

\medskip\par 
In \S \ref{sec: lemmas for lowering the pole order}-\S \ref{sec: finishing the proof} below we handle the computation of the coefficients of 1 and $\bG$ in the linear form associated with $\omega_{F,t}$ under the hypotheses of Theorem \ref{thm: linear forms}, i.e., assuming (i) $F\in \Q[x,y]^\sigma$ and (ii) if $t$ is even (resp. odd), $x^ty^t$ (resp. $x^{t+1}y^{t+1}$) divides $F$. Our approach to compute the coefficients of 1 and $\bG$ under these hypotheses is to reduce the order of the pole to 1, the case for which we already know how to compute the coefficients. The process for the reduction of the order of the pole is elementary, but we will need to be careful about preserving $\sigma$-invariance and regularity in the process. Condition (ii) of Theorem \ref{thm: linear forms} will allow us to maintain nice regularity properties throughout the process. The pole order reduction process will be explicit, in that combined with the explicit calculations of the coefficients of 1 and $\bG$ in the simple pole case, we can explicitly compute the coefficients of 1 and $\bG$ in the linear form coming from every $F$ and $t$ as in Theorem \ref{thm: linear forms}.
\medskip\par 
Finally, in \S \ref{sec: computation in the general integrable case} we will see that the computation of the coefficients of 1 and $\bG$ in the linear form associated with an arbitrary $\omega_{F,t}$ as in Theorem \ref{thm: linear forms, ver2} can be reduced to the case considered in Theorem \ref{thm: linear forms}.

\subsection{Formulas for lowering the order of the pole}\label{sec: lemmas for lowering the pole order}
In this subsection we gather the basic formulas needed for the pole order reduction process. We start by introducing a convenient notation.
\begin{notation}\label{notations for operators}
For any $F\in\Q(x,y)$, set
\[
D_1F:=(\partial/\partial x)(F/x)
\]
and
\[
D_2F:=(\partial/\partial y)(F/y).
\]
\end{notation}
For referencing, we record some useful properties of the operators $D_1$ and $D_2$.
\begin{lemma}\label{lem: D_1 and D_1} \ 
\begin{enumerate}[label=$(\alph*)$] 
\item\label{a} The operators $D_1$ and $D_2$ commute.
\item If $f\in\Q(x,y)$ and $t$ is a non-negative integer, then 
\begin{align*}
D_1(x^ty^tf) & = ((t-1)x^{t-2}f+x^{t-1}f_x) \, y^t\\
D_2(x^ty^tf) & = ((t-1)y^{t-2}f+y^{t-1}f_y) \, x^t
\end{align*}
and
\[
D_2D_1 (x^ty^tf) = x^{t-2}y^{t-2}\bigm((t-1)^2f+(t-1)(xf_x+yf_y)+xyf_{xy}\bigm),
\]
where $f_x$, $f_y$, $f_{xy}=(f_x)_y$ refer to the partial derivatives. In particular, if $F\in \Q[x,y]$ and $x^ty^t\mid F$ (in $\Q[x,y]$) for an integer $t\geq 1$ (resp. $t\geq 2$), then $(D_1+D_2)F$ and $D_2D_1F$ are polynomials (resp. polynomials divisible by $x^{t-2}y^{t-2}$).
\item\label{c} We have
\begin{align*}
D_2D_1 \, \phi_{m,n} & = (m-1)(n-1)\,\phi_{m-2,n-2} \\
D_2D_1 \, \psi_{m,n} & = (m-1)(n-1)\,\psi_{m-2,n-2} \\
(D_1+D_2) \, \phi_{m,n} & = (m-1)\,\phi_{m-2,n}+(n-1)\,\phi_{m, n-2} \\
(D_1+D_2) \, \psi_{m,n} & = (m-1)\,\psi_{m-2,n}+(n-1)\,\psi_{m,n-2},
\end{align*}
where as before, $\phi_{m,n}=x^my^n+x^ny^m$ and $\psi_{m,n}=x^my^n-x^ny^m$ for $m,n\in \Z$.
\item\label{d} If $F\in \Q[x,y]$ is fixed by $\sigma$ and divisible by $xy$, then $D_2D_1F$ and $(D_1+D_2)F$ are polynomials fixed by $\sigma$.
\end{enumerate}
\end{lemma}
\begin{proof}
Parts \ref{a}-\ref{c} are obtained by trivial computations. Part \ref{d} follows from \ref{c} in view of Lemma \ref{lem: basis for sigma-invariants}.
\end{proof}
\medskip\par
Next, we record the formulas for pole reduction. 
\begin{lemma}\label{lem: pole reduction formulas}
Let $t$ be an integer $\geq1$ and $F\in\Q(x,y)$. One has the following equalities of meromorphic differential forms on $\wX$ and $X$.
\begin{enumerate}[label=$(\alph*)$]
\item\label{a'} $\displaystyle{\frac{F}{g^{t+1}}dxdy = \frac{1}{2t} \ d\left(\frac{F}{xg^t}dy\right) - \frac{1}{2t}\frac{D_1F}{g^t}dxdy}$
\item\label{b'} $\displaystyle{\frac{F}{g^{t+1}}dxdy = \frac{-1}{2t} \ d\left(\frac{F}{yg^t}dx\right) - \frac{1}{2t}\frac{D_2F}{g^t}dxdy}$
\item\label{c'} Let $t\geq 2$. Then
\[\frac{F}{g^{t+1}}dxdy = \frac{1}{2t}d\left(\frac{F}{xg^t}dy\right) + \frac{1}{4t(t-1)}d\left(\frac{D_1F}{yg^{t-1}}dx\right) + \frac{1}{4t(t-1)}\frac{D_2D_1F}{g^{t-1}}dxdy.\]
\end{enumerate}
\end{lemma}
\begin{proof}
Part \ref{a'} and \ref{b'} are verified by direct computations. Part \ref{c'} is obtained by successive use of the first two parts.
\end{proof}

\subsection{Calculation of the coefficients of $1,\bG$ in Theorem \ref{thm: linear forms} -- the case of even $t$}\label{sec: proof of linear forms thm, t=even}
Suppose $t$ is an integer $\geq 2$. Let $F\in \Q[x,y]$ be fixed by $\sigma$ and divisible by $x^ty^t$. Our goal is to compute the coefficients of 1 and $\bG$ in the integral of $\omega_{F,t}$ over $\D$. We use the formula of Lemma \ref{lem: pole reduction formulas} ($c$). By Lemma \ref{lem: D_1 and D_1} ($d,b$), $D_2D_1F$ is a $\sigma$-invariant element of $\Q[x,y]$ divisible by $x^{t-2}y^{t-2}$, hence $\omega_{D_2D_1F,t-2}$ is integrable and its integral over $\D$ is in $\langle 1,\bG\rangle_\Q$. Consider the first two terms on the right-hand side of the formula of Lemma \ref{lem: pole reduction formulas} (c), that is, the following sum:
\begin{equation}\label{eq10}
\frac{1}{2t}d\left(\frac{F}{xg^t}dy\right) + \frac{1}{4t(t-1)}d\left(\frac{D_1F}{yg^{t-1}}dx\right) = \frac{1}{4t(t-1)} \ d\left(\frac{2(t-1)x^{-1}Fdy+y^{-1}(D_1F)gdx}{g^t}\right).
\end{equation}
It is worth noting that each of the two summands on the left individually may have a pole along $E_Q\setminus \widetilde{A}$ (this can be seen in coordinate charts). However, their sum is regular on $\widetilde{X}\setminus \widetilde{A}$ as the other two terms in the formula of Lemma \ref{lem: pole reduction formulas} ($c$) are such. In view of the formula of Lemma \ref{lem: pole reduction formulas} ($c$), the following lemma (see also its proof) reduces the computation of the coefficients of $1$ and $\bG$ in the integral of $\omega_{F,t}$ over $\D$ to the computation of these coefficients in the integral of $\omega_{D_2D_1F,t-2}$ over $\D$. When $t$ is even, we may iterate this process until the calculation of the coefficients of 1 and $\bG$ is reduced to the simple pole case. In view of our computations in \S \ref{sec: lin forms II}, this makes Theorem \ref{thm: linear forms} constructive for all even $t$.

\begin{lemma}\label{lem: reduction by 2 leftover term}
Let $t$ be an integer $\geq 2$ and $F\in \Q[x,y]$ be divisible by $x^ty^t$. Then the integral of \eqref{eq10} over $\Delta$ is a rational number.
\end{lemma}
\begin{proof}
We shall compute the integral in question. For simplicity, set
\[
\theta:= \frac{2(t-1)x^{-1}Fdy+y^{-1}(D_1F)gdx}{g^t}
\]
and write $F=x^ty^tf$ (note that both $F/x$ and $(D_1F)/y$ are polynomials). A computation in charts shows that $\theta$ is regular on $\widetilde{X}\setminus \widetilde{A}$. For instance, on the chart $U_3=\{(y, u=\beta_2/\alpha_2)\}$ of the proof of Lemma \ref{deRham} (also recall the coordinates \eqref{p1eq12}), a direct computation shows
\begin{align}\label{p2eq16}
\begin{split}
   \theta \ = \ &\frac{(-1)^{t-1}u^{t-2}y^{t-1}\bigm((t-1)f+u(y-1)f_x\bigm)}{\widetilde{g}^{t-1}} \, du \\
   + \ & \frac{(-1)^{t-1}u^{t-1}y^{t-1}\bigm((t-1)(u^2-1)f +f_x u\widetilde{g}\bigm)}{\widetilde{g}^t}\, dy,
\end{split}
\end{align}
where $\widetilde{g}=1+y-u^2(1-y)$ is the defining polynomial of $\wA$ (and $f=f(x,y)=f(u(y-1),y)$).

Since $\theta$ is holomorphic in $(\wX^\an\setminus \wA^\an) \supset \wD$, we have
\begin{align}\label{p2eq14}
\int_\Delta d\theta \ =  \int_\wD d\theta \ = \int_{\partial \wD} \theta \  \stackrel{(\ast)}{=} \int_{\partial \D} \theta \ + \int_{E_P^\an \, \cap \,\wD} \theta \ + \int_{E_Q^\an \, \cap \,\wD} \theta,
\end{align}
where in $(\ast)$ we divided the boundary of $\wD$ into three pieces, namely, the part in $E_P^\an$, the part in $E_Q^\an$, and the rest, and used functoriality to replace the latter by the integral over $\partial \D$ (the parts $E_P^\an \, \cap \, \wD$ and $E_Q^\an \, \cap \, \wD$ are to be oriented appropriately as in Figure \ref{fig: partial wD} in \S \ref{subsubsec: integrals of 1st type}).

We start by calculating the integral of $\theta$ over $\partial\D$, which is absolutely convergent on every line segment in $\partial\D$ by regularity of $\theta$ on $\wD$. Thanks to the hypotheses that $x^ty^t\mid F$ and $t\geq 2$, we have $x\mid x^{-1}F$ and $y\mid y^{-1}(D_1F)$, so that the 1-forms $x^{-1}Fdy$ and $y^{-1}(D_1F)gdx$ both vanish on the horizontal and vertical parts of $\partial \Delta$. Denoting the straight line segment from $P$ to $Q$ by $[P,Q]$, thus
\[
\int_{\partial \Delta} \theta \ = \int_{[P,Q]} \theta.
\]
Parametrizing $[P,Q]$ by $(1-y,y)$ for $y$ from 0 to 1, on $[P,Q]$ we have
\[
y^{-1}(D_1F)gdx = -2\bigm((t-1)f+(1-y)f_x\bigm)y^{t}(1-y)^{t-1}dy
\]
and
\[
2(t-1)x^{-1}Fdy = 2(t-1)(1-y)^{t-1}y^tfdy,
\]
so that
\[
\theta= \frac{2(t-1)x^{-1}Fdy + y^{-1}(D_1F)gdx}{g^t} = \frac{-2f_x \, y^t(1-y)^tdy}{2^ty^t(1-y)^t} = \, \frac{-f_x\, dy}{2^{t-1}}
\]
on $[P,Q]$. Thus 
\begin{equation}\label{p2eq15}
\int_{\partial\D} \theta \ = \frac{-1}{2^{t-1}} \ \int\limits_0^1 f_x(1-y,y)\,dy.
\end{equation}
We now turn our attention to the integral of $\theta$ over $E_Q^\an\cap \wD$. Recall from the proof of Lemma \ref{Betti} that in our coordinates \eqref{p1eq12}, we have
\[E_Q^\an\cap \wD \ = \{Q\}\times\{[1:-1]\}\times\{[1:m] \mid -1\leq m \leq 0\}.\]
This is contained in the chart $U_3=\{(y, u=\beta_2/\alpha_2)\}$ of the proof of Lemma \ref{deRham}, on which $\theta$ is given by \eqref{p2eq16} above. Since $y=1$ on $E_Q$, we get
\[
\theta \,\big\vert\,_{U_3\,\cap \, E_Q} \, = \frac{(-1)^{t-1}\,(t-1)\,f(Q)}{2^{t-1}}\, u^{t-2} \,du.
\]
We thus get
\begin{equation}\label{p2eq17}
    \int_{E_Q^\an\,\cap\,\wD} \theta \ = \int_{-1}^0 \frac{(-1)^{t-1}\,(t-1)\,f(Q)}{2^{t-1}}\, u^{t-2} \,du \,= \ - \frac{f(Q)}{2^{t-1}}. 
\end{equation}
On the other hand, a computation in chart $U_2$ of the proof of Lemma \ref{deRham} shows that $\theta$ vanishes on $E_P\setminus \wA$. Putting this together with \eqref{p2eq14}, \eqref{p2eq15} and \eqref{p2eq17}, the claim follows since $f\in\Q[x,y]$.
\end{proof}

\begin{rem}
Suppose the polynomial $F=x^ty^tf$ of the last lemma has coefficients in $\Z$. Setting $N=\deg(F)$, then $f_x(1-y,y)\in\Q[y]$ has coefficients in $\Z$ and degree at most $N-2t-1$. It follows from the above argument that the lowest denominator of the integral of \eqref{eq10} over $\Delta$ divides
\[
2^{t+1}t(t-1)L_{N-2t},
\]
where as before, $L_m$ means the least common multiple of $1,2,\ldots, m$.
\end{rem}

\subsection{Calculation of the coefficients of $1,\bG$ in Theorem \ref{thm: linear forms} -- the case of odd $t$}\label{sec: finishing the proof}

In this subsection we will make the remaining case of Theorem \ref{thm: linear forms} constructive. Suppose $t\geq 1$ is an odd integer and $F\in \Q[x,y]$ is $\sigma$-invariant and divisible by $x^{t+1}y^{t+1}$. Averaging the first two formulas of Lemma \ref{lem: pole reduction formulas}, we have
\begin{equation}\label{eq11}
\frac{F}{g^{t+1}}dxdy = \frac{1}{4t} \left\{ d\left( \frac{F\,dy}{xg^t}\right) - d \left(\frac{F\,dx}{yg^t}\right) - \frac{(D_1+D_2)F}{g^t}dxdy \right\}.
\end{equation}
By Lemma \ref{lem: D_1 and D_1}, $(D_1+D_2)F$ is a $\sigma$-invariant polynomial divisible by $x^{t-1}y^{t-1}$. Thus the integral of $\omega_{(D_1+D_2)F,t-1}$ over $\D$ is in $\langle 1,\bG\rangle_\Q$, and this integral can be written as an explicit linear combination of $1$ and $\bG$ using \S \ref{sec: lin forms II} and the procedure of \S \ref{sec: proof of linear forms thm, t=even}. Hence, the following lemma (see also its proof) will complete the task of making Theorem \ref{thm: linear forms} constructive for all $t$.

\begin{lemma}
Let $t$ be an integer $\geq 1$ and $F\in \Q[x,y]$ be divisible by $x^{t+1}y^{t+1}$. Then each of the integrals
\[
\int_\Delta d\left( \frac{F\,dy}{xg^t}\right) \ \ \ \ \text{and} \ \ \ \ \int_\Delta d \left(\frac{F\,dx}{yg^t}\right)
\]
is convergent and its value is a rational number.
\end{lemma}
\begin{proof}
We only compute the first integral as the other integral is computed similarly. A computation in charts shows that since $x^{t+1}y^{t+1}$ divides $F$, the differential 1-form $\frac{Fdy}{xg^t}$ is regular on $\widetilde{X}\setminus\widetilde{A}$. (Note that this form may have a pole along $E_Q\setminus \widetilde{A}$ if we only assume that $x^t$ divides $F$.) One easily sees that $\frac{Fdy}{xg^t}$ vanishes on $E_P\setminus \wA$ and $E_Q\setminus \wA$. It follows that we have an equality of absolutely convergent integrals
\[
\int_\Delta d\left( \frac{F\,dy}{xg^t}\right) =  \int_{\partial \Delta} \frac{F\,dy}{xg^t}.
\]
The condition $x^{t+1}\mid F$ guarantees that $F/x$ is a multiple of $x$ and hence the 1-form $\frac{Fdy}{xg^t}$ vanishes on the vertical line segment in $\partial \Delta$. Of course, it also vanishes on the horizontal line segment, hence
\[
\int_{\partial \Delta} \frac{F\,dy}{xg^t} = \int_{[P,Q]} \frac{F\,dy}{xg^t}
\]
(where $[P,Q]$ is again the straight line segment from $P$ to $Q$). Writing $F=x^{t+1}y^{t+1}f$ and parametrizing $[P,Q]$ by $(1-y,y)$ for $y$ between 0 and 1, the last integral reads
\[
\int_0^1 \frac{(1-y)^{t+1}y^{t+1}f(1-y,y)\,dy}{2^t(1-y)^{t+1}y^t} = \frac{1}{2^t} \int_0^1 yf(1-y,y)\,dy,
\]
which is a rational number.
\end{proof}

\subsection{Calculation of the coefficients for arbitrary integrable differential forms}\label{sec: computation in the general integrable case}
In this final subsection, we show that the computation of the coefficients of 1 and $\bG$ for an arbitrary $\omega_{F,t}$ as in Theorem \ref{thm: linear forms, ver2} can be simplified to the computations of these coefficients for differential forms satisfying the hypotheses of Theorem \ref{thm: linear forms}. The new ingredient is the following lemma.
\begin{lem}\label{lem: thmv1 implies thmv2}
Let $F\in \Q[x,y]$ and $t\in \Z_{\geq 0}$. Suppose that $F$ is fixed by both $\sigma$ and $\tau$, and that $\omega_{F,t}$ is integrable. Then, modulo the space of regular 2-forms on $X$, the differential form $\omega_{F,t}$ on $X$ is in the $\Q$-span of the differential forms $\omega_{x^my^m, s}$, where $m,s$ are non-negative integers, $m$ is even, and $s\leq m$. (Note that by Lemma \ref{lem: 1st char. of order}, for a differential form $\omega_{x^my^m, s}$ with $m,s\in\Z_{\geq 0}$, the condition $s\leq m$ is equivalent to integrability.)
\end{lem}
\begin{proof}
The space $\Q[x,y]^{\sigma, \tau}$ is equal to the space of all symmetric polynomials in $x^2$ and $y^2$. It follows that $F$ is a polynomial in $x^2y^2$ and $x^2+y^2$ with rational coefficients, or equivalently, $F$ is a polynomial in $x^2y^2$ and $g$ with rational coefficients. Start by writing
\begin{equation}\label{p2eq20}
F = \sum\limits_{\substack{m,r\geq 0\\ m\,=\,\text{even}}} e_{m,r} \, (xy)^mg^r,
\end{equation}
where each $e_{m,r}$ is in $\Q$. Then
\[
\omega_{F,t} \,= \sum\limits_{\substack{t\geq r\geq 0\\ m \,=\,\text{even}\,\geq 0}} e_{m,r} \, \omega_{x^my^m,\, t-r} \ + \ \text{a regular 2-form on $X$}.
\]
It thus suffices to show that whenever $e_{m,r}$ is nonzero, we have $m\geq t-r$. Suppose $F$ is not zero and $\nu$ is the smallest value of $m+r$ such that $e_{m,r}\neq 0$ for some $m,r$. To compute the order of $F$ along $E_P$, we write
\[
F(x,u(x-1)) \, = \sum\limits_{\substack{r\geq 0\\ m \,= \,\text{even}\,\geq 0}} (-1)^re_{m,r} \, x^mu^m(x-1)^{m+r} \,\widetilde{g}^r,
\]
where $g(x,u(x-1))=(1-x)\widetilde{g}$ with $\widetilde{g}=1+x-u^2(1-x)$. Separating the pairs $(m,r)$ with $m+r=\nu$ from the rest, we can write
\begin{align*}
F(x,u(x-1)) \, = & \, (x-1)^\nu \sum\limits_{\nu\geq \,m \,= \, \text{even}\,\geq 0} (-1)^{\nu-m}e_{m,\nu-m} \, x^mu^m \widetilde{g}^{\nu-m} \\
& + (x-1)^{\nu+1} H(x,u)
\end{align*}
for some $H(x,u)\in \Q[x,u]$. In view of Lemma \ref{lem: 1st char. of order}, from the last equation it follows that $\ord_{E_P}F=\nu$, as the polynomial
\[
\sum\limits_{\nu\geq \,m \,= \, \text{even}\,\geq 0} (-1)^{\nu-m}e_{m,\nu-m} \, x^mu^m \widetilde{g}^{\nu-m} \ \in \Q[x,u]
\]
does not vanish in $\Q[u]$ at $x=1$. Since $\omega_{F,t}$ is integrable, it follows that $\nu\geq t$, as desired.
\end{proof}

Now suppose $F\in\Q[x,y]$ is fixed by $\sigma$ and $\omega_{F,t}$ is integrable. Our goal is to find the coefficients of 1 and $\bG$ in the integral of $\omega_{F,t}$ over $\D$. We first argue that we may replace $F$ by its $(\sigma,\tau)$-invariant component: Write $F=F_++F_-$ according to the decomposition \eqref{p2eq2}, where $F_+\in \Q[x,y]^{\sigma, \tau}$ and $F_-\in \Q[x,y]^{\substack{\sigma=1 \\ \tau=-1}}$. We claim that $\omega_{F_+, t}$ and $\omega_{F_-, t}$ are both integrable. Indeed, the involution $\tau$ acts on $X$ and hence on the space of meromorphic 2-forms on $X$, and the decomposition $\omega_{F,t}=\omega_{F_+,t}+\omega_{F_-,t}$ is the decomposition of $\omega_{F,t}$ according to the eigenspaces of $\tau$, with $\tau$ acting on $\omega_{F_-,t}$ as 1. Our claim follows from the fact that the space of integrable 2-forms on $X$ is invariant under $\tau$. (The $\tau$-invariance of this space follows from Lemma \ref{lem: integrability condition}, since $\tau$ extends to an automorphism of $\wX$ that stabilizes $\wX\setminus\wA$, see Remark \ref{rem: automorphism tau}.) Since the integral of $\omega_{F_-, t}$ over $\D$ is zero (see Lemma \ref{lem: sigma=1, tau=-1 F} or Remark \ref{rem: interpretation of tau-anti invariance lem}), we have
\[
\int_\D \omega_{F,t} = \int_\D \omega_{F_+,t}.
\]
Thus to compute the coefficients of 1 and $\bG$, we may replace $F$ by $F_+$.

Next, by Lemma \ref{lem: thmv1 implies thmv2} and its proof, we can write $\omega_{F_+, t}$ explicitly as the sum of a regular form $\omega'$ on $X$ and a $\Q$-linear combination of differentials of the form $\omega_{x^my^m, s}$ with $m,s\in\Z_{\geq 0}$, $m$ even, and $s\leq m$. The integral of $\omega'$ (which is a polynomial in $\Q[x,y]$ times $dxdy$) on $\D$ is a rational number. Each of the forms $\omega_{x^my^m, s}$ with $m,s\in\Z_{\geq 0}$, $m$ even, and $s\leq m$ satisfies the hypotheses of Theorem \ref{thm: linear forms}. Hence we know how to compute the coefficients of 1 and $\bG$ in the integral of each of these forms over $\D$.

\begin{rem}
For the sake of completeness, here we record another version of Theorem \ref{thm: linear forms, ver2}: for every $t\in\Z_{\geq 0}$ and $F\in \Q[x^2,y^2]$, if $\omega_{F,t}$ is integrable, then the integral of $\omega_{F,t}$ over $\D$ is in $\langle 1,\bG\rangle_{\Q}$. The equivalence of this statement with Theorem \ref{thm: linear forms, ver2} follows from the following three facts: (i) For any $F\in \Q[x,y]$ with integrable $\omega_{F,t}$, writing $F=F_++F_-$ according to the decomposition \eqref{p2eq22}, both $\omega_{F_+,t}$ and $\omega_{F_-,t}$ are integrable; (ii) the integral of $\omega_{F_-,t}$ over $\D$ is zero; and (iii) both $\Q[x^2,y^2]$ and $\Q[x,y]^\sigma$ contain $\Q[x,y]^{\sigma,\tau}$, and the images of $\Q[x^2,y^2]$ and $\Q[x,y]^\sigma$ under the eigenspace projection map $\Q[x,y]\rightarrow \Q[x,y]^\tau$ are both $\Q[x,y]^{\sigma,\tau}$.
\end{rem}

\begin{rem}
While we arrived at the statement and proof of Theorem \ref{thm: linear forms, ver2} through motivic considerations, the computations of \S \ref{sec: lin forms II} and \S \ref{sec: lin forms III} can also be used to give a constructive and elementary (albeit perhaps unmotivated and seemingly out of nowhere, especially in the simple pole case) proof of Theorem \ref{thm: linear forms, ver2}, with no reference to the underlying motivic context. Indeed, the computations of \S \ref{sec: lin forms III} reduce the problem to the simple pole case. In the simple pole case, given $F\in \Q[x,y]^\sigma$ one first {\it defines} $b_F$ by the formula of Proposition \ref{prop: computation of b_F}, and then defines $a_F$ by the formula
\[
a_F = \int_\Delta \omega_F-b_F\,\omega_0,
\]
so that $\int_\Delta \omega_F=a_F+b_F\,\bG$ (where $\omega_F=\omega_{F,0}$, as before). Then one proceeds to show that $a_F$ is rational and explicitly computable through the computations of \S \ref{sec: coefficient of eta, log case}.
\end{rem}

\section*{Acknowledgments}
We are grateful to Michel Waldschmidt, Francis Brown, and Daniel Bertrand for insightful discussions and correspondences. We are also grateful to the anonymous referees for their careful reading of the paper and many comments and suggestions, which helped us improve the exposition of the paper. 

This article was written while Nemoto was visiting the University of Winnipeg. He would like to thank the Department of Mathematics and Statistics for their hospitality. Nemoto was supported by Waseda University Grant for Early Career Researchers (Project number: 2025E-041). Eskandari and Nemoto were partially supported by the NSERC Discovery Grant ``Arithmetic and Geometry of Mixed Motives" (RGPIN-2024-04076).

\end{document}